\documentclass[a4paper,10pt]{article}
\usepackage[english]{babel}
\usepackage[utf8]{inputenc}

\usepackage{amsmath}
\usepackage{amsfonts}
\usepackage{amssymb}
\usepackage{amsthm}
\usepackage{graphicx}
\usepackage{pdfpages}

\newcommand{\mytext}[1]{ \: \textrm{#1} \: }
\newcommand{\mysetdescr}[2]{\left\{ { #1 \: \left| \: #2 \right. } \right\} }

\newcommand{\myN}{\mathbb{N}}
\newcommand{\myNk}[1]{\underline {#1}}
\newcommand{\myR}{\mathbb{R}}
\newcommand\mytimes{{\times}}

\newcommand{\strH}{\sqsubseteq_\H}

\newcommand{\strS}{\sqsubseteq_\S}

\def\A{{\cal A}}

\def\H{{\cal H}}
\def\I{{\cal I}}
\def\J{{\cal J}}
\def\L{{\cal L}}
\def\M{{\cal M}}
\def\P{{\cal P}}

\def\S{{\cal S}}
\def\Z{{\cal Z}}

\newcommand{\mf}[1]{\mathfrak{ #1 }}
\newcommand{\fB}{\mf{B}}
\newcommand{\fC}{\mf{C}}
\newcommand{\fD}{\mf{D}}

\newcommand{\fo}{\mf{o}}
\newcommand{\fP}{\mf{P}}
\newcommand{\fR}{\mf{R}}
\newcommand{\fT}{\mf{T}}
\newcommand{\fU}{\mf{U}}
\newcommand{\fz}{\mf{z}}

\newcommand{\fDr}{\fD_r}
\newcommand{\fTa}{\fT_a}

\newcommand{\noth}{{- \! \! \! \! h}}

\newcommand{\fTahn}[1]{ \fTa^{({#1})} }
\newcommand{\fTan}{ \fTahn{n} }
\newcommand{\fTahnA}[1]{ \fTa^{< {#1} >} }
\newcommand{\fTanA}{ \fTahnA{n} }
\newcommand{\fChn}[1]{ \fC^{[{#1}]} }
\newcommand{\fCn}{ \fChn{n} }

\newcommand{\Nin}{N^{in}}
\newcommand{\Nout}{N^{out}}

\newcommand{\iotaxi}[1]{\iota_{#1}}
\newcommand{\iotaMxi}[2]{\iotaxi{#2, #1}}

\def\BP{\begin{proof}}
\def\EP{\end{proof}}

\DeclareMathOperator{\id}{id}

\DeclareMathOperator{\Tr}{Tr}
\DeclareMathOperator{\Rd}{Rd}

\begin{document}

\theoremstyle{plain}
\newtheorem{condition}{Condition}
\newtheorem{theorem}{Theorem}
\newtheorem{definition}{Definition}
\newtheorem{corollary}{Corollary}
\newtheorem{lemma}{Lemma}
\newtheorem{proposition}{Proposition}

\title{\bf Interdependencies of less-equal-relations between partial Lov\'{a}sz-vectors of digraphs}
\author{\sc Frank a Campo}
\date{\small Seilerwall 33, D 41747 Viersen, Germany\\
{\sf acampo.frank@gmail.com}}

\maketitle

\begin{abstract}
\noindent For digraphs $G$ and $H$, let $\H(G,H)$ be the set of all homomorphisms from $G$ to $H$, and let $\S(G,H)$ be the subset of those homomorphisms mapping all proper arcs in $G$ to proper arcs in $H$. From an earlier investigation we know that for certain digraphs $R$ and $S$, the relation ``$\# \S(G,R) \leq \# \S(G,S)$ for all $G \in \fD'$'' implies ``$\# \H(G,R) \leq \# \H(G,S)$ for all $G \in \fD'$'', where $\fD'$ is a subclass of digraphs. Now we ask for the inverse: For which digraphs $R, S$ and which subclasses $\fD'$ of digraphs does ``$\# \H(G,R) \leq \# \H(G,S)$ for all $G \in \fD'$'' imply ``$\# \S(G,R) \leq \# \S(G,S)$ for all $G \in \fD'$''? We prove this implication for three combinations of digraph classes. In particular, the relations are equivalent for all flat posets $R, S$ with respect to all flat posets $G$.
\newline

\end{abstract}

\section{Introduction} \label{sec_introduction}

The number of homomorphisms between directed graphs (digraphs) may carry important information about structure. Graph parameters which can be expressed as numbers of homomorphisms into weighted graphs have been characterized by Freedman et al.\ \cite{Freedman_etal_2007} in 2007. In different fields of graph theory \cite{Hell_Nesetril_2004,Schroeder_2016}, the still open reconstruction conjecture asks if two objects with at least four vertices are isomorphic if all numbers of embeddings of certain subgraphs into them are equal. 

A classical result is the Theorem of Lov\'{a}sz \cite{Lovasz_1967} from 1967. It states that numbers of homomorphisms distinguish non-isomorphic general ``relational structures''; with $\H(G,H)$ defined as the set of homomorphisms from a digraph $G$ to a digraph $H$, the following specifications are relevant for our purpose:
\newpage

\begin{theorem}[Lov\'{a}sz \cite{Lovasz_1967}] \label{theo_Lovasz}
Let $\fD'$ be the class of finite digraphs or the class of finite posets. Then, for $R, S \in \fD'$
\begin{align*}
R & \simeq S \\
\Leftrightarrow \; \; \; \# \H(G,R) & = \# \H(G,S) \; \mytext{for every } G \in \fD'.
\end{align*}
For the class of finite posets, the equivalence holds also if we replace the homomorphism sets by the sets of strict order homomorphisms.
\end{theorem}

A short and simple proof for digraphs is contained in \cite{Hell_Nesetril_2004}; with minor modification, it works for posets, too.

The infinite vector $\L_\H(H) \equiv ( \# \H(G,H) )_{G \in \fD'}$ is called the {\em Lov\'{a}sz-vector of $H$}. In the last decades, topics related to it have found interest \cite{Borgs_etal_2006,Lovasz_2006,Freedman_etal_2007,
Borgs_etal_2008,Lovasz_Szegedy_2008,Schrijver_2009,
Borgs_etal_2012,Cai_Govorov_2020} in connection with vertex and edge weights. For undirected graphs, Dvo\v{r}\'{a}k \cite{Dvorak_2010} investigated in 2010 proper sub-classes $\fU'$ of undirected graphs for which the partial Lov\'{a}sz-vector $( \# \H(G,H) )_{G \in \fU'}$ still distinguishes graphs; the distinguishing power of the vector $( \# \H(G,H) )_{H \in \fU'}$ has been investigated by Chaudhuri and Vardi \cite{Chaudhuri_Vardi_1993} in 1993 and Fisk \cite{Fisk_1995} in 1995. The significance of the Lov\'{a}sz-vector for other relations than isomorphism has been studied by Dell et al.\ \cite{Dell_etal_2018} in 2018 and Atserias et al.\  \cite{Atserias_etal_2021} in 2021.

\begin{figure} 
\begin{center}
\includegraphics[trim = 40 760 235 35, clip]{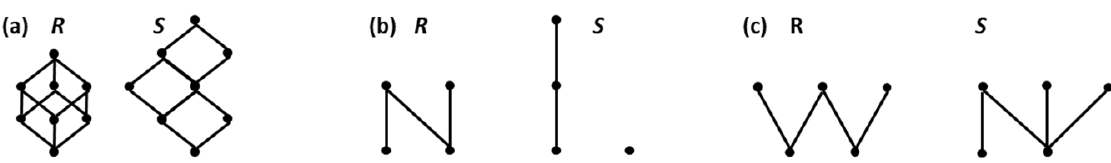}
\caption{\label{fig_Intro} Hasse-diagrams of three pairs of posets $R$ and $S$ with $\# \H(P,R) \leq \# \H(P,S)$ for every finite poset $P$.}
\end{center}
\end{figure}

This paper continues the work of the author about the pointwise less-equal-relation between partial Lov\'{a}sz-vectors of digraphs:
\begin{quote}
{\em Given a class $\fD'$ of finite digraphs, what is it in the structure of finite digraphs $R$ and $S$ that enforces}
\end{quote}
\begin{equation} \label{fragestellung}
\# \H(G,R) \leq \# \H(G,S) \; \mytext{\em for every } G \in \fD'  \mytext{\em ?}
\end{equation}

The starting point of the work was the pair of partial ordered sets (posets) $R$ and $S$ shown in Figure \ref{fig_Intro}(a). The author \cite[Theorem 5]{aCampo_2018} has proven that for these posets $\# \H(P,R) \leq \# \H(P,S)$ for every finite poset $P$. Additional non-trivial examples are shown in the Figures \ref{fig_Intro}(b)-(c); more are contained in \cite{aCampo_toappear_0,aCampo_toappear_1}.

In order to have a compact notation, we define (using the abbreviation ``w.r.t.'' for ``with respect to''):

\begin{definition} \label{def_sqsubseteq}
Let $R, S$ be digraphs and let $\fD'$ be a class of digraphs. We write
\begin{align*}
R \strH S \mytext{ w.r.t. } \fD' & \; \mytext{ iff } \;
\# \H(G,R) \leq \# \H(G,S) \mytext{ for all } G \in \fD', \\
R \strS S  \mytext{ w.r.t. } \fD' & \; \mytext{ iff } \;
\# \S(G,R) \leq \# \S(G,S) \mytext{ for all } G \in \fD',
\end{align*}
where $\S(G,T)$ is the set of {\em strict} homomorphisms from $G$ to $T$, i.e., the set of homomorphisms mapping all proper arcs in $G$ to proper arcs in $T$. We call $\fD'$ the {\em domain} of the relations $\strH$ and $\strS$.
\end{definition}

In \cite{aCampo_toappear_0,aCampo_toappear_1}, the author has started a systematic examination of the relations $R \strH S$ and $R \strS S$ for several domains. The subject of \cite{aCampo_toappear_0} were sufficient criteria for $R \strH S $, whereas \cite{aCampo_toappear_1} focused on the theoretical aspects of $R \strS S$. One of the main results of \cite{aCampo_toappear_0} was the following theorem in which $\fD, \fTa$, $\fP$ and $\fP^*$ denote the classes of all digraphs, of the digraphs with antisymmetrical transitive hull, of the posets, and of the posets with loops removed:

\begin{theorem}[\cite{aCampo_toappear_0}, Theorem 2] \label{theo_GschemeOnStrict}
Let $R \in \fD$ and $\fD' \subseteq \fD$. Then, for all $S \in \fD$, the implication
\begin{align} \label{RsqGS_then_RsqS}
R \strS S \quad \Rightarrow \quad R \strH S
\end{align}
with domain $\fD'$ on both sides holds if
\begin{align*}
R & \in \fD' = \fD, \\
R & \in \fTa \subseteq \fD' \subseteq \fD, \\
\mytext{or} \quad R & \in \fD' \mytext{ with } \fD' = \fP \mytext{ or } \fD' = \fP^*.
\end{align*}
\end{theorem}
With $\L_\S(H) \equiv ( \# \S(G,H) )_{G \in \fD'}$, this Theorem  describes conditions under which a point wise less-equal-relation between the Lov\'{a}sz-vectors $\L_\S(R)$ and $\L_\S(S)$ implies a point wise less-equal-relation between the Lov\'{a}sz-vectors $\L_\H(R)$ and $\L_\H(S)$. In this article, we ask for the inverse implication:
\begin{quote}
{\em For which classes $\fD_0, \fD', \fR$ of digraphs do we have
\begin{align} \label{zielimplikation}
R \strH S \mytext{ w.r.t. } \fD_0 & \quad \Rightarrow \quad R \strS S \mytext{ w.r.t. } \fD'
\end{align}
for all $R, S \in \fR$? }
\end{quote}
We thus ask for triplets of digraph classes: on the one hand, we have a class $\fR$ from which we take $R$ and $S$, and on the other hand, we have the domains $\fD_0$ and $\fD'$ of the relations $\strH$ and $\strS$. It is the appropriate combination of these classes which makes the implication \eqref{zielimplikation} working, and in the Theorems \ref{theo_R_R_mytimes}, \ref{theo_fTan}, and \ref{theo_fTanA}, we prove it for three combinations.

Our approach is the following. For $G \in \fD'$, we construct a sequence $G_\nu$ of digraphs in $\fD'$ in such a way that $\# \H(G_\nu, T) $ is an exponential function with exponent $\nu$ for every $T \in \fR$. For an appropriate choice of $\fD'$ and $\fR$, the leading term of the exponential function is dominated by the number of strict homomorphisms from $G$ to $T$. Because $R \strH S$ enforces $\# \H(G_\nu,R) \leq \# \H(G_\nu,S)$ for all $\nu$, the relation $R \strH S$ yields $\# \S(G,R) \leq \# \S(G,S)$ for $G \in \fD'$, hence $R \strS S$ with respect to $\fD'$.

After recalling fundamental terms of graph theory in Section \ref{sec_preparation}, we forge our tools in Section \ref{sec_arcWeights} by means of {\em arc weights} which are instances of the edge weights mentioned above: mappings from $A(G^*)$ to $\myN_0$, where $A(G^*)$ is the set of proper arcs in the digraph $G$. For an arc weight $\alpha$ of $G$, we construct a sequence $G(\alpha)_\nu$ of digraphs for which $\# \H( G(\alpha)_\nu,R)$ is an exponential function with exponent $\nu$ for every reflexive digraph $R$. We show in Theorem \ref{theo_selExists}, how {\em selecting arc weights} can be constructed for certain ``extremal'' homomorphisms $\xi \in \H(G,R)$; for a $\xi$-selecting arc weight $\alpha$ of $G$, the leading term of $\# \H(G(\alpha)_\nu,R)$ is determined by the cardinality of a homomorphism set $ \J_{G,R}(\xi)$ closely related to $\xi$. Based on this result, it is shown in Theorem \ref{theo_IGR_IGS} that for reflexive digraphs $R, S \in \fD' \in \{ \fD, \fTa, \fP \}$, the relation $R \strH S$ with respect to $\fD'$ implies $ \# \J_{G,R}(\xi) \leq \# \J_{G,S}(\zeta) $ for related extremal homomorphisms $\xi \in \H(G,R)$, $\zeta \in \H(G,S)$.

In the implications of type \eqref{zielimplikation} proven in this article, all involved digraphs belong to $\fTa$. In Section \ref{sec_fTa}, we look deeper into this class. We discuss the properties of paths in digraphs $G \in \fTa$, and we recall the definition of the {\em transitive reduction} of $G$. By means of the transitive reduction of paths, we define a class $\fR \subset \fTa$ of reflexive digraphs, and we see in Lemma \ref{lemma_classR}, that for $R \in \fR$ and $G \in \fTa$ with $\S(G,R) \not= \emptyset$, the existence of a suitable extremal strict homomorphism is guaranteed.

Theorem \ref{theo_IGR_IGS} and Lemma \ref{lemma_classR} are combined and applied in Section \ref{sec_implications}, resulting in three implications of type \eqref{zielimplikation}. Theorem \ref{theo_R_R_mytimes} in Section \ref{subsec_R_R_mytimes} concerns digraphs $R, S \in \fR$ with ``sparse'' arc sets. In Theorem \ref{theo_fTan} in Section \ref{subsec_fTan}, it is shown that for digraphs $R, S \in \fR$ of equal height, the relation $R \strH S$ with respect to $\fTa$ (with respect to $\fP)$ implies $R \strS S$ with respect to a certain subclass of $\fTa$ (subclass of $\fP$) characterized by maximal paths. In particular, the relations $R \strH S$ and $R \strS S$ with domain $\fP$ or $\fTa$ are equivalent for flat posets $R$ and $S$. By weakening the requirements on $G$ and tightening those on $R$ in Section \ref{subsec_fTanA}, we show in Theorem \ref{theo_phiG} that for this new combination of digraph properties, the relation between $\# \S(G,R)$ and $\# \J_{G,R}(\sigma)$ does not depend on $R$ for strict homomorphisms $\sigma$. Using this result, we present in Theorem \ref{theo_fTanA} a variant of Theorem \ref{theo_fTan} referring to the new combination of digraph classes.

Section \ref{sec_lemma_fTan_invers} contains a proof which has been postponed because the result does not belong to the main line of this article.

\section{Preparation} \label{sec_preparation}

\subsection{Basics and Notation} \label{subsec_notation}

A {\em (finite) directed graph} or {\em digraph} $G$ is an ordered pair $(V(G),A(G))$ in which $V(G)$ is a non-empty, finite set and $A(G) \subseteq V(G) \mytimes V(G)$ is a binary relation on $V(G)$. We write $vw$ for an ordered pair $(v,w) \in V(G) \mytimes V(G)$. The elements of $V(G)$ are called the {\em vertices} of $G$ and the elements of $A(G)$ are called the {\em arcs} of $G$. A digraph $G$ is {\em reflexive},  or {\em symmetric}, or {\em antisymmetric}, etc., iff the relation $A(G)$ has the respective property. A {\em partially ordered set (poset)} is a reflexive, antisymmetric, transitive digraph.

For digraphs $G$ and $H$, we write $H \subseteq G$ and call $H$ a {\em subgraph} of $G$, iff $V(H) \subseteq V(G)$ and $A(H) \subseteq A(G)$. For $\emptyset \not= X \subseteq V(G)$, the {\em digraph $G \vert_X$ induced on $X$} has vertex set $X$ and arc set $A(G) \cap ( X \mytimes X)$.

An arc $vw \in A(G)$ is called {\em proper} iff $v \not= w$; otherwise, it is called a {\em loop}. All possible loops of $G$ are collected in the {\em diagonal} $\Delta_G \equiv \mysetdescr{(v,v)}{v \in V(G)}$. $G^* \equiv (V(G), A(G) \setminus \Delta_G)$ is the digraph $G$ {\em with loops removed}. Vertices $v, w \in V(G)$ are {\em adjacent} iff $vw \in A(G)$ or $wv \in A(G)$. The {\em neighborhood} $N_G(v)$ of $v \in V(G)$ is the set of all $w \in V(G)$ adjacent to $v$. Furthermore,
\begin{align*}
\Nin_G(v) & \; \equiv \; \mysetdescr{ w \in N_G(v) }{ wv \in A(G) }, \\
\Nout_G(v) & \; \equiv \; \mysetdescr{ w \in N_G(v) }{ vw \in A(G) }.
\end{align*}
A vertex $v \in V(G)$ is {\em isolated} iff $N(v) \subseteq \{ v \}$.

For $G \in \fD$ and every $vw \in A(G)$, the {\em interval $[ v, w ]_G$} is defined as
\begin{equation*}
[ v, w ]_G \quad \equiv \quad \Nout_G(v) \cap \Nin_G(w),
\end{equation*}
and the cardinality of $[ v, w ]_G$ is denoted by
\begin{equation*}
\iota(v,w)_G \quad \equiv \quad \# [ v, w ]_G.
\end{equation*}
We regard $\iota$ as mapping from $A(G)$ to $\myN_0$. If $G$ is reflexive, then $v, w \in [ v, w ]_G$ for all $vw \in A(G)$.

Vertices $v, w \in V(G)$ are {\em connected in $G$} iff $v = w$ or - in the case of $v \not= w$ - there exists a sequence of vertices $z_0, \ldots , z_I$ in $V(G)$ with $v = z_0$, $w = z_I$, $I \in \myN$, and $z_{i-1}, z_i$ adjacent for all $1 \leq i \leq I$. We say that the {\em line} $z_0, \ldots , z_I$ in $V(G)$ {\em connects} $v$ and $w$. A non-empty subset $X \subseteq V(G)$ is {\em connected} iff all its vertices are connected in $G \vert_X $, and it is called a {\em connectivity component} of $G$, iff it is connected and $X' = X$ holds for every connected subset $X' \subseteq V(G)$ with $X \subseteq X'$.

A sequence $W = v_0, \ldots , v_I$ of vertices of $G$ with $I \in \myN_0$ is called a {\em (directed) walk} iff $v_{i-1}v_i \in A(G)$ for all $1 \leq i \leq I$. The {\em length of $W$} is $\ell(W) \equiv I$, and we say that $W$ {\em starts in $v_0$} and {\em ends in $v_I$}, or that $W$ is a walk {\em from $v_0$ to $v_I$}. If all vertices in $W$ are distinct, $W$ is called a {\em path in $G$}. $W$ is {\em closed} iff $I \geq 1$ and $v_0 = v_I$. A digraph is {\em acyclic} iff it does not contain a closed walk. We define $G^i$, $i \in \myN$, as the graph with vertex set $V(G)$ and $vw \in A(G^i)$ iff there exists a walk in $G$ of length $i$, starting in $v$ and ending in $w$.

For a digraph $G$, the {\em transitive hull} $\Tr(G)$ of $G$ is the digraph with vertex set $V(G)$ and the (set-theoretically) smallest transitive arc set containing $A(G)$:
\begin{align*}
A(\Tr(G)) \quad & \equiv \quad \bigcap \mysetdescr{ T \subseteq V(G) \mytimes V(G) }{ T \mytext{ transitive with } A(G) \subseteq T }.
\end{align*}
$vw \in A(\Tr(G))$ is equivalent to the existence of a walk in $G$ from $v$ to $w$:
\begin{align*}
A(\Tr(G)) \quad & = \quad \bigcup_{i=1}^\infty A( G^i ).
\end{align*}
$\Tr(G)$ contains all loops of $G$.

We need symbols for several classes of digraphs. $\fD$ is the class of all digraphs with finite non-empty vertex set, $\fDr \subset \fD$ is the class of the reflexive digraphs, and $\fP \subset \fDr$ is the class of the finite posets.

The following class will play an important role in Sections \ref{sec_fTa} and \ref{sec_implications}:
\begin{align*}
\fTa & \; \equiv \; \mysetdescr{ G \in \fD }{ G^* \mytext{is acyclic} }.
\end{align*}
Equivalently, $\fTa$ can be characterized as the class of digraphs with antisymmetric transitive hull (that is the reason for the choice of the symbol $\fTa$) or as the class of subgraphs of posets. In particular, all digraphs in $\fTa$ are antisymmetric, and for $G \in \fTa$, also the transitive hull of $G$ belongs to $\fTa$. Important sub-classes contained in $\fTa$ are $\fP$ and (trivially) the class of acyclic digraphs. In Lemma \ref{lemma_path_in_fTa} in Section \ref{subsec_Paths}, we will characterize the class $\fTa$ by paths.

From set theory, we use the following notation:
\begin{align*}
\myNk{0} & \equiv  \emptyset, \\
\myNk{n} & \equiv  \{ 1, \ldots, n \} \mytext{for every} n \in \myN, \\
\myNk{n}_0 & \equiv \myNk{n} \cup \{ 0 \} \mytext{for every} n \in \myN_0.
\end{align*}

On the vertex set $\myNk{n}_0$, $n \in \myN_0$, we define the poset $C_n$ by setting $A(C_n) \equiv \mysetdescr{ ij \in \myNk{n}_0 \mytimes \myNk{n}_0 }{ i \leq j }$. It is known as {\em chain of length $n$} in order theory, and $C_n^*$ is called {\em transitive tournament} in graph theory.

Let $I \in \myN$, let $a_1, \ldots , a_I$ be positive real numbers, and let $0 < x_1 < \ldots < x_I$ be pairwise different positive real numbers. We call a mapping 
\begin{align*}
f : \myN_0 & \rightarrow \myR \\
\nu & \mapsto \sum_{i=1}^I a_i \cdot x_i^\nu
\end{align*}
an {\em exponential function} (with exponent $\nu$), and we call $( a_I, x_I )$ its {\em leading term}. For two exponential functions $f$ and $g$ with leading terms $( a, x )$ and $( b, y )$, respectively, we have $f(\nu) > g(\nu)$ for all sufficiently large values of $\nu$ if $x > y$ or if $ x = y $ and $a > b$.

$\A(X,Y)$ is the set of mappings from $X$ to $Y$. For $f \in \A(X,Y)$ and $X' \subseteq X$, we write $f \vert_{X'}$ for the pre-restriction of $f$ to $X'$, and for any set $Z$ with $f(X) \subseteq Z$, we write $f \vert^{Z}$ for the post-restriction of $f$ to $Z$.

For any set $X$, $\id_X : X \rightarrow X$, $x \mapsto x$, is the identity mapping of $X$, the mapping $\fo_X : X \rightarrow \myN_0$ is the {\em zero-mapping} with $\fo_X(x) = 0$ for all $x \in X$, and $\fz_X : X \rightarrow \myN_0$ maps every $x \in X$ to 2.

Finally, we use the {\em Cartesian product}. Let $\I$ be a non-empty set, and let $M_i$ be a non-empty set for every $i \in \I$. Then the Cartesian product of the sets $M_i, i \in \I$, is defined as
\begin{eqnarray*}
\prod_{i \in \I} M_i & \; \equiv \; & 
\mysetdescr{ f \in \A \big( \I, \bigcup_{i \in \I} M_i \big)}{ f(i) \in M_i \mytext{for all} i \in \I }.
\end{eqnarray*}

\subsection{Homomorphisms} \label{subsec_homomorphisms}

Given digraphs $G$ and $H$, we call a mapping $\xi : V(G) \rightarrow V(H)$ a {\em homomorphism from $G$ to $H$} if $\xi(v) \xi(w) \in A(H)$ for all $vw \in A(G)$. For such a mapping, we write $\xi : G \rightarrow H$; we collect the homomorphisms in the set
\begin{align*}
\H(G,H) & \; \equiv \; \mysetdescr{ \xi : V(G) \rightarrow V(H) }{ \xi \mytext{ is a homomorphism} }.
\end{align*}
Isomorphism is indicated by ``$\simeq$''.

Let $G, G' \in \fD$ with $G \subseteq G'$. We summarize closely related homomorphisms in $\H(G',H)$ into sets by defining for all $H \in \fD$ and all $\xi \in \H(G,H)$
\begin{equation*}
[ \xi ]_{G'} \; \equiv \; \mysetdescr{ \xi' \in \H(G',H) }{ \xi = \xi' \vert_{V(G)}  }.
\end{equation*}
For every $\xi' \in \H(G',H)$, we have $\xi' \vert_{V(G)} \in \H(G,H)$ and $ \xi' \in \left[ \xi' \vert_{V(G)} \right]_{G'}$. Furthermore, $[ \xi ]_{G'} \cap [ \zeta ]_{G'} = \emptyset$ for $\xi, \zeta \in \H(G,H)$ with $\xi \not= \zeta$. Therefore, 
$\mysetdescr{ [ \xi ]_{G'} }{ \xi \in \H(G,H) }$ is a partition of $\H( G', H )$ iff $[\xi]_{G'} \not= \emptyset$ for all $\xi \in \H(G,H)$, and in this case $\H(G,H)$ is a representation system of the partition. 

Every homomorphism $\xi : G \rightarrow H$ maps loops in $G$ to loops in $H$, but proper arcs of $G$ can be mapped to both, loops and proper arcs of $H$. We call a homomorphism from $G$ to $H$ {\em strict} iff it maps all proper arcs of $G$ to proper arcs of $H$.
\begin{align*}
\S(G,H) & \; \equiv \; \H(G,H) \cap \H(G^*,H^*)
\end{align*}
is the set of strict homomorphisms from $G$ to $H$. The set $\H(G^*,H^*) \setminus \H(G,H)$ contains all homomorphisms from $G^*$ to $H^*$ which map a vertex belonging to a loop in $G$ to a vertex of $H$ not belonging to a loop. A mapping $\xi : V(G) \rightarrow V(H)$ is thus a strict homomorphism, iff it maps loops in $G$ to loops in $H$ and proper arcs of $G$ to proper arcs of $H$. If $H$ is reflexive, then $\S(G,H) = \H(G^*,H^*)$; for posets $G$ and $H$, the set $\S(G,H) = \H(G^*,H^*)$ is the set of strict order homomorphisms from $G$ to $H$.

For $G, H \in \fD$ and $\xi \in \H(G,H)$, we use the notation
\begin{align*}
\iotaxi{\xi}(v,w) & \quad \equiv \quad \iota( \xi(v),\xi(w))_H \quad \mytext{for all } vw \in A(G), \\
\mu_\xi(G) & \quad \equiv \quad \sum_{vw \in A(G^*)} \iotaxi{\xi}(v,w).
\end{align*}
We regard also $\iotaxi{\xi}$ as mapping from $A(G)$ to $\myN_0$. If $H$ is a reflexive digraph from $\fTa$, then, for all $\xi \in \H(G,H)$,
\begin{equation} \label{iocaChar_strict}
\xi \in \S(G,H) \quad \Leftrightarrow \quad
\iotaxi{\xi}(v,w) \geq 2 \; \; \mytext{for all} \; vw \in A(G^*).
\end{equation}

In what follows, our focus is frequently on the pre-restriction of $\iotaxi{\xi}$ to subsets of $A(G^*)$. In order to unburden the notation, we define for all $G, H \in \fD$, $B \subseteq A(G^*)$, and all $\xi \in \H(G,H)$
\begin{equation*}
\iotaMxi{B}{\xi} \quad \equiv \quad \left( \iotaxi{\xi} \right) \vert_{B}.
\end{equation*}

\section{Arc weights} \label{sec_arcWeights}

\subsection{Definition and basic properties} \label{subsec_def_arcWeight}

\begin{figure}
\begin{center}
\includegraphics[trim = 40 740 295 55, clip]{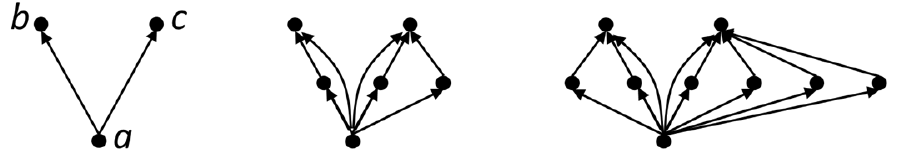}
\caption{\label{fig_Def_Ganu} A digraph $G$ and the digraphs $G(\alpha)_1$ and $G(\alpha)_2$ resulting for the arc weight $\alpha$ with $\alpha(a,b) = 1$ and $\alpha(a,c) = 2$.}
\end{center}
\end{figure}

Our main instrument in this study are {\em arc weights} and sequences of digraphs constructed by means of them:

\begin{definition} \label{def_arcweight}
For every $G \in \fD$, we call a mapping $\alpha : A(G^*) \rightarrow \myN_0$ an {\em arc weight} of $G$, and we define
\begin{align*}
D(\alpha) \; & \equiv \; \mysetdescr{ vw \in A(G^*) }{ \alpha(v,w) > 0  }.
\end{align*}

Given an arc weight $\alpha$ of a digraph $G \in \fD$, we construct the digraphs $G(\alpha)_\nu$, $\nu \in \myN_0$, in the following way: Let $X_\nu(v,w)$, $vw \in A(G^*)$, be pairwise disjoint sets of cardinality $\nu \cdot \alpha(v,w)$, all disjoint from $V(G)$. We define
\begin{align*}
V( G(\alpha)_\nu ) & \; \equiv \; V(G) \; \cup \bigcup_{vw \in D(\alpha)} X_\nu(v,w), \\
A( G(\alpha)_\nu ) & \; \equiv \; A(G) \; \cup \bigcup_{vw \in D(\alpha)} \left( \{ v \} \mytimes X_\nu(v,w) \right) \cup \left( X_\nu(v,w) \mytimes \{ w \} \right).
\end{align*}
If $G$ is a poset, we denote with $G'(\alpha)_\nu$ the transitive hull of $G(\alpha)_\nu$ with loops added for all vertices contained in the sets $X_\nu(v,w)$, $vw \in D(\alpha)$.
\end{definition}

The digraph $G(\alpha)_\nu$ is thus constructed by clamping a set of $\nu \cdot \alpha(v,w)$ singletons between the vertices $v$ and $w$ with $vw \in D(\alpha)$. An example is shown in Figure \ref{fig_Def_Ganu}. The construction method is a variant of the replacement operation described in \cite[Section 4.4]{Hell_Nesetril_2004}. 

We have $G(\alpha)_0 = G$, and for every $\nu \in \myN_0$, we have $[v, w]_{G(\alpha)_\nu} = [v, w]_G \cup X_\nu(v,w)$ for all $vw \in A(G^*)$ (with $X_\nu(v,w) = \emptyset$ for $vw \notin D(\alpha)$). Later on, we need

\begin{lemma} \label{lemma_Gan_in_Ta}
Let $\alpha$ be an arc weight of $G \in \fTa$. Then $G(\alpha)_\nu \in \fTa$ for all $\nu \in \myN_0$. If $G$ is a poset, then $G'(\alpha)_\nu$ is a poset, too.
\end{lemma}
\BP Let $W = v_0, \ldots , v_I$ be a walk in $G(\alpha)_\nu^*$. If $v_i \notin V(G)$ for an $i \in \myNk{I-1}$, then $v_i \in X_\nu(v_{i-1},v_{i+1})$ with $v_{i-1} v_{i+1} \in A(G^*) \subseteq A(G(\alpha)_\nu^*)$. By skipping all vertices $v_i \notin V(G)$ with $i \in \myNk{I-1}$, we get a walk $W' = w_0, \ldots , w_J$ in $G(\alpha)_\nu^*$ with $w_0 = v_0$, $w_J = v_I$, and $w_j \in V(G)$ for all $j \in \myNk{J-1}$.

Assume that $W$ is a closed walk. Then also $W'$ is a closed walk. $w_0 = w_J \in V(G)$ is not possible due to $G \in \fTa$. Assume $w_0 = w_J \in X_\nu(a,b)$ for $ab \in D(\alpha)$. Then $w_1 = b$, $w_{J-1} = a$, and $a, w_1, \ldots , w_{J-1}$ is a closed walk in $G^*$ which is impossible.

Assume $G \in \fP \subset \fTa$. We have already seen $G(\alpha)_\nu \in \fTa$, and therefore the transitive hull of $G(\alpha)_\nu$ is antisymmetric. It contains a loop for every $v \in V(G)$, and by adding loops for all vertices contained in the sets $X_\nu(v,w)$, $vw \in D(\alpha)$, it becomes a poset.

\EP

Let $H \in \fD$, $\xi \in \H(G,H)$, and let $\zeta \in [ \xi ]_{G(\alpha)_\nu}$, $\nu \geq 1$. The homomorphism $\zeta$ coincides with $\xi$ on $G$. In consequence, it maps the set $X_\nu(v,w)$ to the interval $[ \xi(v), \xi(w) ]_H$ for all $vw \in D(\alpha)$, and $\zeta$ is uniquely determined by the family of mappings $f_{vw} \equiv \zeta \vert_{X_\nu(v,w)} \in \A( X_\nu(v,w), [ \xi(v), \xi(w) ]_H)$, $vw \in D(\alpha)$. On the other hand, if $H$ is reflexive, then every
\begin{equation*}
f \; \in \; \prod_{vw \in D(\alpha)} \A \left( X_\nu(v,w), [ \xi(v), \xi(w) ]_H \right) \quad \mytext{\em (Cartesian product)}
\end{equation*}
defines a homomorphism $\zeta_f \in [ \xi ]_{G(\alpha)_\nu}$ via
\begin{align*}
\zeta_f(a) & \equiv \begin{cases}
\xi(a) & \mytext{if }  a \in V(G), \\
f_{vw}(a) & \mytext{if } a \in X_\nu(v,w).
\end{cases}
\end{align*}
This assignment between the Cartesian product and $[ \xi ]_{G(\alpha)_\nu}$ is bijective, and we conclude that for $H \in \fDr$, we have
\begin{equation}
\label{card_xiGanu}
\# [ \xi ]_{G(\alpha)_\nu} \quad = \quad \pi_\alpha(\xi)^\nu
\end{equation}
where
\begin{equation*}
\pi_\alpha(\xi) \quad \equiv \quad \prod_{vw \in A(G^*)} \iotaxi{\xi}(v,w)^{\alpha(v,w)}.
\end{equation*}

We have $G = G(\alpha)_0$ and $G \subseteq G(\alpha)_\nu$ for every $\nu \in \myN_0$, and in consequence, the set
$\mysetdescr{ [ \xi ]_{G(\alpha)_\nu} }{ \xi \in \H(G,H) }$ is a partition of $\H( G(\alpha)_\nu, H )$ with representation system $\H(G,H)$, hence
\begin{equation} \label{HGanuH_summe}
\# \H( G(\alpha)_\nu, H ) \quad = \quad \sum_{\xi \in \H(G,H)} \# [ \xi ]_{G(\alpha)_\nu}.
\end{equation}
In particular, $\# \H(G(\alpha)_\nu, H)$ is an exponential function with exponent $\nu$.

Also for a poset $G$ and $\xi \in \H(G,H)$, every $\zeta \in [ \xi ]_{G'(\alpha)_\nu}$ maps $X_\nu(v,w)$ to $[ \xi(v),\xi(w)]_H$. In consequence, Equation \eqref{card_xiGanu} remains valid for $H \in \fDr$ if we exchange $G(\alpha)_\nu$ against $G'(\alpha)_\nu$, and $\# \H(G'(\alpha)_\nu, H)$ is an exponential function with exponent $\nu$, too.

Because isolated vertices in $G$ do not belong to any proper arc of $G$, they do not affect the value of $\pi_\alpha(\xi)$ and do not contribute to the size of $[ \xi ]_{G(\alpha)_\nu}$. They can be neglected in what follows. It should be mentioned that the concepts and results presented in this article are still valid (but become trivial) if $G$ consists of isolated vertices only: in this case, there exists only the trivial arc weight $\fo_\emptyset = ( \emptyset, \emptyset, \myN_0 )$ (written as a mapping), and for all $H \in \fDr$ and all $\xi \in \H(G,H)$, we have $\iotaxi{\xi} = ( \emptyset, \emptyset, \myN )$ and $\pi_{\fo_\emptyset}(\xi) = 1$; furtheremore, $\H(G,H) = \S(G,H)$

The following lemma shows how an arc weight of a digraph $G$ can be constructed by referring to arc weights of subgraphs of $G$:

\begin{lemma} \label{lemma_wProdukt}
Let $G \in \fD$, let $\L$ be a set of subgraphs of $G$, and let $\alpha_L$ be an arc weight of $L$ for every $L \in \L$. We define for every $vw \in A(G^*)$
\begin{equation*}
\beta(v,w) \quad \equiv \quad \sum_{\stackrel{L \in \L}{vw \in A(L^*)}} \alpha_L(v,w).
\end{equation*}
Then $ \beta$ is an arc weight of $G$ with $D(\beta) = \cup_{L \in \L} D( \alpha_L)$ and
\begin{align}
\label{formel_beta}
\pi_\beta(\theta) = &
\prod_{L \in \L} \pi_{\alpha_L}( \theta \vert_{L} )
\end{align}
for every $\theta \in \H(G,H)$, $H \in \fDr$.
\end{lemma}
\BP We define for all $L \in \L$ and all $vw \in A(G^*)$
\begin{eqnarray*}
\alpha'_L(v,w)& \equiv & 
\begin{cases}
\alpha_L(v,w), & \mytext{if} \; vw \in A(L^*), \\
0 & \mytext{otherwise.}
\end{cases}
\end{eqnarray*}
Then $\alpha'_L$ is an arc weight of $G$ for all $L \in \L$, and for all $H \in \fDr$, $\theta \in \H(G, H)$, we have
\begin{align*}
\pi_\beta(\theta) = &
\prod_{vw \in A(G^*)} \iotaxi{\theta}(v,w)^{\beta(v,w)} \\
= & \;
\prod_{vw \in A(G^*)} \prod_{L \in \L} \iotaxi{\theta}(v,w)^{\alpha'_L(v,w)} \\
= & \;
\prod_{L \in \L} \prod_{vw \in A(G^*)} \iotaxi{\theta}(v,w)^{\alpha'_L(v,w)} \\
= & \;
\prod_{L \in \L} \prod_{vw \in A(L^*)} \iotaxi{\theta}(v,w)^{\alpha_L(v,w)} \\
= & \;
\prod_{L \in \L} \pi_{\alpha_L}( \theta \vert_{L} ).
\end{align*}

\EP

\subsection{Selecting arc weights} \label{subsec_sel_arcWeights}

\begin{definition} \label{def_selArcW}
Given $G \in \fD$, $H \in \fDr$, and $\zeta \in \H(G,H)$, we call an arc weight $\alpha$ of $G$ {\em $\zeta$-selecting within $\H(G,H)$} iff, for all $\xi \in \H(G,H)$,
\begin{equation} \label{ineq_selective}
\pi_\alpha(\xi) \quad \leq \quad \pi_\alpha(\zeta)
\end{equation}
with equality iff $\iotaMxi{D(\alpha)}{\xi} = \iotaMxi{D(\alpha)}{\zeta}$.
\end{definition}

The zero-mapping $\fo_{A(G^*)}$ is for every $H \in \fDr$ and every $\zeta \in \H(G,H)$ a $\zeta$-selecting arc weight of $G$ within $\H(G,H)$.

The motivation for introducing selecting arc weights becomes visible in the next corollary which results directly from \eqref{HGanuH_summe} and \eqref{card_xiGanu}:
\begin{corollary} \label{coro_Jxi}
For $G \in \fD$, $H \in \fDr$, and $\zeta \in \H(G,H)$, let $\alpha$ be a $\zeta$-selecting arc weight of $G$ within $\H(G,H)$. Then $\# \H(G(\alpha)_\nu,H)$ is an exponential function with exponent $\nu$ and leading term
\begin{equation*}
\left( i_\zeta, \pi_\alpha(\zeta) \right)
\end{equation*}
where $i_\zeta$ is the number of homomorphisms $\xi \in \H(G,H)$ with $\iotaMxi{D(\alpha)}{\xi} = \iotaMxi{D(\alpha)}{\zeta}$. \end{corollary} 

The next lemma shows, that the combination of arc weights described in Lemma \ref{lemma_wProdukt} results in a selecting arc weight if all combined arc weights are selecting.

\begin{lemma} \label{lemma_selecting_1}
Let $G \in \fD$, $H \in \fDr$, $\zeta \in \H(G,H)$, and let $\L$ be a set of subgraphs of $G$. For every $L \in \L$, let $\alpha_L$ be a $\zeta \vert_{L}$-selecting arc weight of $L$ within $\H(L,H)$. With $\beta$ defined as in Lemma \ref{lemma_wProdukt}, $\beta$ is a $\zeta$-selecting arc weight of $G$ within $\H(G,H)$.
\end{lemma}
\BP Because of Lemma \ref{lemma_wProdukt}, $\beta$ is an arc weight of $G$ with $D(\beta) = \cup_{L \in \L} D( \alpha_L)$. For every $\xi \in \H(G,H)$,
\begin{equation*}
\pi_\beta(\xi)
\; \stackrel{\eqref{formel_beta}}{=} \; 
\prod_{L \in \L} \pi_{\alpha_L}( \xi \vert_{L} )
\; \stackrel{\eqref{ineq_selective}}{\leq} \; 
\prod_{L \in \L} \pi_{\alpha_L}( \zeta \vert_{L} )
\; \stackrel{\eqref{formel_beta}}{=} \; 
\pi_\beta(\zeta),
\end{equation*}
with equality iff $\iotaMxi{D(\alpha_L)}{\xi \vert_{L}} = \iotaMxi{D(\alpha_L)}{\zeta \vert_{L}}$ for all $L \in \L$. We have to show that this is equivalent to $\iotaMxi{D(\beta)}{\xi}  = \iotaMxi{D(\beta)}{\zeta}$.

For every $\theta \in \H(G,H)$, we have for every $L \in \L$
\begin{equation} \label{eq_ithetaDL_ithetaD}
\iotaMxi{D(\alpha_L)}{\theta \vert_{L}} \; = \;
\left( \iotaxi{\theta \vert_{L}} \right) \vert_{D(\alpha_L)} \; = \;
\left( \iotaxi{\theta} \right) \vert_{D(\alpha_L)} \; = \;
\iotaMxi{D(\alpha_L)}{\theta },
\end{equation}
and $\iotaMxi{D(\beta)}{\xi} = \iotaMxi{D(\beta)}{\zeta}$ implies $\iotaMxi{D(\alpha_L)}{\xi \vert_{L}} = \iotaMxi{D(\alpha_L)}{\zeta \vert_{L}}$ for all $L \in \L$.

Now assume $\iotaMxi{D(\alpha_L)}{\xi \vert_{L}} = \iotaMxi{D(\alpha_L)}{\zeta \vert_{L}}$ for all $L \in \L$. For $vw \in D(\beta)$, there exists an $L \in \L$ with $vw \in D(\alpha_L)$, hence
\begin{align*}
\iotaMxi{D(\beta)}{\xi }(v, w) \; = \;
\iotaMxi{D(\alpha_L)}{\xi }(v, w) 
\; \stackrel{\eqref{eq_ithetaDL_ithetaD}}{=} \;
\iotaMxi{D(\alpha_L)}{\zeta }(v, w) \; = \
\iotaMxi{D(\beta)}{\zeta }(v, w).
\end{align*}

\EP

Our key in constructing selecting arc weights is {\em Gibb's inequality} from information theory:
\begin{lemma}[{\cite[Lemma 1]{Topsoe_1974}}] \label{lemma_Topsoe}
Let $x_1, \ldots , x_n$ and $y_1, \ldots , y_n$ be positive numbers with $\sum_{i=1}^n x_i \leq \sum_{i=1}^n y_i$. Then
\begin{equation*}
\prod_{i=1}^n {x_i}^{y_i} \quad \leq \quad \prod_{i=1}^n {y_i}^{y_i} 
\end{equation*}
with equality iff $x_i = y_i$ for every $i \in \myNk{n}$.
\end{lemma}

In order to work in a concise manner with this lemma and its consequences, we need:
\begin{definition} \label{def_hatmu_usw}
Let $G, L \in \fD$, $H \in \fDr$, and let $\L$ be a set of subgraphs of $G$. We define
\begin{align*}
\hat{\mu}(L,H) & \; \equiv \; \max \mysetdescr{ \mu_\xi(L) }{ \xi \in \H(L,H) }, \\
\M(L,H) & \; \equiv \; \mysetdescr{ \xi \in \H(L,H) }{ \mu_\xi(L) = \hat{\mu}(L,H) }, \\
\M^\L(G,H) & \; \equiv \; \mysetdescr{ \xi \in \H(G,H) }{ \xi \vert_{V(L)} \in \M(L,H) \mytext{ for all } L \in \L  }, \\
\I^\L(G,H) & \; \equiv \; \mysetdescr{ \iotaMxi{\cup_{L \in \L} A(L^*)}{\xi} }{ \xi \in \M^\L(G,H) }.
\end{align*}
\end{definition}
Observe that $\hat{\mu}$ is invariant with respect to isomorphism: $L \simeq L'$ and $H \simeq H'$ yields $\hat{\mu}(L,H) = \hat{\mu}(L',H')$.

\begin{corollary} \label{coro_sel_MLH}
Let $L \in \fD$, $H \in \fDr$, and $\zeta \in \M(L, H)$. Then $\iotaMxi{A(L^*)}{\zeta}$ is a $\zeta$-selecting arc weight of $L$ within $\H(L,H)$.
\end{corollary}
\BP Let $\xi \in \H(L,H)$, $B \equiv A(L^*)$. Then
\begin{align*}
\sum_{vw \in B} \iota_\xi(v,w) \; = \; \mu_\xi(L) \; \leq \; \hat{\mu}(L,H) \; = \; \sum_{vw \in B} \iota_\zeta(v,w).
\end{align*}
In the case of $B = \emptyset$, everything is trivial, and in the case of $B \not= \emptyset$, we conclude with Lemma \ref{lemma_Topsoe}
\begin{align*}
\pi_{\iotaMxi{B}{\zeta}}(\xi) & = \prod_{vw \in B} {\iotaxi{\xi}(v,w)}^{\iotaxi{\zeta}(v,w)} \\
& \leq 
\prod_{vw \in B} {\iotaxi{\zeta}(v,w)}^{\iotaxi{\zeta}(v,w)} \; = \; \pi_{\iotaMxi{B}{\zeta}}(\zeta),
\end{align*}
with equality iff $\iotaMxi{B}{\xi} = \iotaMxi{B}{\zeta}$.

\EP

\begin{theorem} \label{theo_selExists}
Let $G \in \fD$, $H \in \fDr$, and let $\L$ be a set of subgraphs of $G$. For every $\zeta \in \M^\L(G,H)$, we get a $\zeta$-selecting arc weight $\gamma_\zeta^\L$ of $G$ within $\H(G,H)$ by defining for all $vw \in A(G^*)$
\begin{align*}
\gamma_\zeta^\L(v,w) \; & \equiv \; \iotaxi{\zeta}(v,w) \cdot \# \mysetdescr{ L \in \L }{ vw \in A(L^*) }.
\end{align*}
\end{theorem}
\BP According to Corollary \ref{coro_sel_MLH}, the mapping $ \alpha_L \equiv \iotaMxi{A(L^*)}{\zeta}$ is for every $L \in \L$ a $\zeta \vert_L$-selecting arc weight of $L$ within $\H(L,H)$. For $vw \in A(G^*)$,
\begin{align*}
\gamma_\zeta^\L(v,w) \; & = \;
\sum_{\stackrel{L \in \L}{(v, w) \in A(L^*)}} \iotaxi{\zeta}(v,w)  \; = \; 
\sum_{\stackrel{L \in \L}{(v, w) \in D(\alpha_L)}} \alpha_L(v,w).
\end{align*}
Now apply Lemma \ref{lemma_selecting_1}.

\EP

\begin{definition} \label{def_Jxi}
Let $G \in \fD$, $H \in \fDr$, and let $\L$ be a set of subgraphs of $G$. For $\xi \in \H(G,H)$ and $H' \in \fDr$, we define
\begin{align*}
\J_{G,H'}^\L(\xi) & \; \equiv \; \mysetdescr{ \zeta \in \H(G,H') }{ \iotaMxi{A(L^*)}{\zeta} = \iotaMxi{A(L^*)}{\xi} \mytext{ for all } L \in \L}.
\end{align*}
\end{definition}
For all $\zeta \in \J_{G,H'}^\L(\xi)$, $\xi \in \H(G,H)$, we have  $\J_{G,H'}^\L(\zeta) = \J_{G,H'}^\L(\xi)$. Furthermore, $\J_{G,H'}^\emptyset(\xi) = \H(G,H')$ for all $\xi \in \H(G,H)$.

\begin{theorem} \label{theo_IGR_IGS}
Let $R \in \fDr \cap \fD'$ with $\fD' \in \{ \fP, \fTa, \fD \}$, and let $S \in \fDr$ with $R \strH S$ with respect to $\fD'$. Furthermore, let $G \in \fD'$ and let $\L$ be a set of subgraphs of $G$, $B \equiv \cup_{L \in \L} \; A(L^*)$. Then, for all $\xi \in \M^\L(G,R)$,
\begin{align*}
\iotaMxi{B}{\xi} \in \I^\L(G,S) \quad \Rightarrow \quad 
\# \J_{G,R}^\L(\xi) & \leq \# \J_{G,S}^\L(\xi).
\end{align*}
\end{theorem}

\BP Let $\xi \in \M^\L(G,R)$ with $\iotaMxi{B}{\xi} \in \I^\L(G,S)$. There exists a $\zeta \in \M^\L(G,S)$ with $\iotaMxi{B}{\zeta} = \iotaMxi{B}{\xi}$, hence $\gamma_\xi^\L = \gamma_\zeta^\L$. Due to Theorem \ref{theo_selExists}, $\delta \equiv \gamma_\xi^\L$ is an arc weight of $G$ which is $\xi$-selecting within $\H(G,R)$ and $\zeta$-selecting within $\H(G,S)$.

According to Corollary \ref{coro_Jxi}, $\# \H(G(\delta)_\nu,R)$ is an exponential function with exponent $\nu$ and leading term $\left( \# \J_{G,R}^\L(\xi), \pi_\delta(\xi) \right)$, and $\# \H(G(\delta)_\nu,S)$ is an exponential function with exponent $\nu$ and leading term $\left( \# \J_{G,S}^\L(\zeta), \pi_\delta(\zeta) \right)$. But due to $\iotaMxi{B}{\zeta} = \iotaMxi{B}{\xi}$, we have $\J_{G,S}^\L(\zeta) = \J_{G,S}^\L(\xi)$ and $ \pi_\delta(\zeta) = \pi_\delta(\xi)$. The leading term of $\# \H(G(\delta)_\nu,S)$ is thus $ \left( \# \J_{G,S}^\L(\xi), \pi_\delta(\xi) \right)$. 

In the case of $\fD' = \fD$, the inequality $\J_{G,R}^\L(\xi) \leq \# \J_{G,S}^\L(\xi)$ is now a direct consequence of $R \strH S$. For $\fD' = \fTa$, observe Lemma \ref{lemma_Gan_in_Ta}, and for $\fD' = \fP$, exchange $G(\delta)_\nu$ against $G'(\delta)_\nu$ and use Lemma \ref{lemma_Gan_in_Ta} again.

\EP

\section{Digraphs belonging to $\bf \fTa$} \label{sec_fTa}

In the implications of type \eqref{zielimplikation} proven later on in this article, all involved digraphs belong to $\fTa$. In this section, we look deeper into the properties of paths in such digraphs and we recall the definition of the {\em transitive reduction}. Finally, we define a class $\fR \subset \fTa$ of reflexive digraphs which guarantees the existence of a strict homomorphism with a selecting arc weight (if a strict homomorphism exists at all).

\subsection{Paths} \label{subsec_Paths} \label{subsec_paths}

Digraphs belonging to $\fTa$ have nice properties making it easy to work with them. One is that the {\em concatenation} of two paths is a path again. Let $P = p_0, \ldots , p_I$ and $Q = q_0, \ldots , q_J$ be paths in $G \in \fTa$ with $p_I = q_0$. Then $W \equiv p_0, \ldots , p_I, q_1, \ldots , q_J$ is a walk in $G^*$. $p_i = q_j$ is not possible for $i \in \myNk{I}_0$ and $j \in \myNk{J}$, because $G^*$ is acyclic. Therefore, $W$ is a path in $G$.

Fundamental is:

\begin{lemma} \label{lemma_path_in_fTa}
Let $G \in \fD$ be a digraph. Then $G$ is an element of $\fTa$, iff every path $P = z_0, \ldots , z_I$ in $G$ is uniquely determined by its vertex set ${\mysetdescr{ z_i }{ i \in \myNk{I}_0}}$.
\end{lemma}
\BP Let $G \in \fTa$. For a path $P = z_0, \ldots , z_I$ in $G$, we have to show that the sequence $z_0, \ldots , z_I$ is determined by the set $Z_0 \equiv \mysetdescr{ z_i }{ i \in \myNk{I}_0}$.  Because $G^*$ is acyclic, we have $ z_0 \notin \mysetdescr{ z_i }{ i \in \myNk{I}}$, and we can identify $z_0$ as the only point $v \in Z_0$ with $\Nin_G(v) \cap Z_0 \subseteq \{ v \}$. By proceeding with the set $Z_1 \equiv \mysetdescr{ z_i }{ i \in \myNk{I}}$, we identify  $z_1$, and so on.

If $G^*$ is not acyclic, then there exists a shortest closed walk $W = z_0, \ldots , z_I$ in $G^*$. We have $I \geq 2$, and due to the minimal length of $W$, the sequences $P = z_0, \ldots , z_{I-1}$ and $P' \equiv z_1, \ldots , z_I$ are different paths in $G$ with identical vertex sets.

\EP

According to this lemma, we can identify a path in $G \in \fTa$ with the set of its vertices. In what follows, ``$P \subseteq V(G)$ is a path'' means that there exists a path in $G$ with vertex set $P$. With $P_\bot$ and $P_\top$, we denote the starting point and ending point of the path, respectively, and we notate the sequence belonging to it by indexing its symbol: $P = P_0, \ldots, P_{\ell(P)}$. The largest length of a path in $G$ is called the {\em height of $G$} and notated by $h_G$; digraphs with $h_G \leq 1$ are called {\em flat}.

We call a path $P \subseteq V(G)$ {\em maximal} iff $P = P'$ holds for every path $P' \subseteq V(G)$ with $P \subseteq P'$. A path of length $h_G$ is always maximal.

\begin{corollary} \label{coro_maxPath}
Let $P \subseteq V(G)$ be a maximal path of $G \in \fTa$. For every $i \in \myNk{\ell(P)}$, the two-element sequence $P_{i-1}, P_i$ is the only path in $G$ from $P_{i-1}$ to $P_i$. In particular, if $G$ is additionally reflexive,
\begin{equation} \label{maxPath_iotaFunktion}
\iota(P_{i-1}, P_i)_G \; = \; 2 \quad \mytext{for all } i \in \myNk{\ell(P)},
\end{equation} 
\end{corollary}
\BP Let $Q \subseteq V(G)$ be a path from $P_{i-1}$ to $P_i$ in $G$. Then
\begin{equation*}
P_0, \ldots , P_{i-1}, Q_1, \ldots Q_{\ell(Q)-1} , P_i, \ldots, P_{\ell(P)}
\end{equation*}
is a path in $G^*$, and the maximality of $P$ yields $\ell(Q)  \leq 1$, hence $Q = P_{i-1}, P_i$.  \eqref{maxPath_iotaFunktion} follows, because the sequence $P_{i-1}, v, P_i$ is a path in $G$ for every $v \in [ P_{i-1}, P_i ]_G \setminus \{ P_{i-1}, P_i \}$.

\EP

We define:

\begin{itemize}
\item $\P_G$ is the set of subsets $P \subseteq V(G)$ which are maximal paths in $G$;
\item $\P_G^h$ is the set of the $P \in \P_G$ with $\ell(P) = h_G$;
\item $ V(G)^h \equiv \cup \P_G^h$ is the set of vertices belonging to the paths of length $h_G$;
\item $G^h \equiv G \vert_{V(G)^h}$ is the subgraph of $G$ consisting of the vertices belonging to paths of length $h_G$ in $G$;
\item $ V(G)^\noth \equiv V(G) \setminus V(G)^h$  is the set of vertices not belonging to any path of length $h_G$.
\end{itemize}
In the case of $h_G = 0$, we have $G^h = G$, but in the case of $h_G \geq 1$, the digraph $G^h$ does not contain any isolated vertex.

The following observation will be the key in several proofs:

\begin{lemma} \label{lemma_gFunktion}
For $G \in \fTa$, there exists a mapping $g : V(G)^h \rightarrow \myNk{h_G}_0$ with the following property: For every $P \in \P_G^h$, we have $v = P_{g(v)}$ for all $v \in P$.
\end{lemma}
\BP Let $P, P' \in \P_G^h$ and $v \in P \cap P'$. There exist $i, j \in \myNk{h_P}_0$ with $P_i = v = P'_j$. In the case of $i > j$, the sequence $P_0, \ldots, P_i, P'_{j+1}, \ldots P'_{h_G}$ is a path in $G$ with length greater than $h_G$. In the same way we treat the case $j > i$, and $i = j$ is proven.

\EP

For $G, H \in \fTa$ and $\sigma \in \S(G,H)$, the image $\sigma[P]$ of a path $P$ in $G$ is a path in $H$ with the same length as $P$. In particular, $\S(G,H) \not= \emptyset$ implies $h_G \leq h_H$. Furthermore, for $n \geq h_G$,
\begin{align*}
\lambda_{G,n} : G & \rightarrow C_n, \\ \nonumber
v & \mapsto \max \mysetdescr{ i \in \myNk{n}_0 }{ v = P_i \mytext{ for a path } P \subseteq V(G) }, \\
\mytext{and} \quad \hat{\lambda}_{G,n} : G & \rightarrow C_n, \\
v & \mapsto
\begin{cases}
\lambda_{G,n}(v), & \mytext{if} \Nout_G(v) \not= \emptyset, \\
n, & \mytext{otherwise}
\end{cases}
\end{align*}
are strict homomorphisms. Lemma \ref{lemma_fTan} below shows that the mapping $g$ from Lemma \ref{lemma_gFunktion} is in fact a strict homomorphism from $G^h$ to $\myNk{h_G}_0$, and it is easily seen that $\lambda_{G,h_G}$ is a strict extension of $g$ to $G$.

\subsection{The transitive reduction} \label{subsec_transRed}

The counterpart of the transitive hull of a digraph $G$ is its {\em transitive reduction} $\Rd(G)$. It has been defined by Aho et al.\ \cite{Aho_etal_1972} as a digraph $H$ with vertex set $V(G)$ and the smallest number of arcs resulting in $G \subseteq \Tr(H)$. For $G \notin \fTa$, the arc set of $\Rd(G)$ must not be uniquely determined, and $\Rd(G)$ must not be a subgraph of $G$; a construction method is described in \cite[Theorem 2]{Aho_etal_1972}. But in the case of $G \in \fTa$, the arc set of the transitive reduction of $G$ is according to \cite[Theorem 1]{Aho_etal_1972} given by
\begin{align*}
\bigcap \mysetdescr{ B \subseteq A(G) }
{ G \subseteq \Tr((V(G),B)) },
\end{align*}
and $ \Rd(G)$ is a subgraph of $G$. Alternatively, we get the arc set of $\Rd(G)$ by erasing all arcs $vw \in A(G)$ for which a path of length greater/equal 2 exists in $G^*$, starting in $v$ and ending in $w$. Therefore, the arc set of $\Rd(G)$ can also be written as
\begin{equation} \label{RdG_formel}
A(G) \setminus \bigcup_{i=2}^\infty A((G^*)^i).
\end{equation}
With
\begin{equation*}
G_\mytimes \; \equiv \; \Rd( G )^*,
\end{equation*}
$G_\mytimes$ is a subgraph of $G$. For a poset $G$, $A(G_\mytimes)$ is known as the {\em cover relation} of $G$; it contains exactly the arcs drawn in the Hasse-diagram of $G$. (In fact, the description \eqref{RdG_formel} of the arc set of $G_\mytimes$ for a poset $G$ is widely used and known in order theory independently from Aho et al.\ \cite{Aho_etal_1972}; see, e.g., \cite[p.\ 77]{Erne_1982}.) 

\begin{corollary} \label{coro_pathsInGmytimes}
Let $G \in \fTa$, $P \subseteq V(G)$.
\begin{align*}
P \mytext{path in} G_\mytimes & \quad \Rightarrow \quad P \mytext{path in} G; \\
P \mytext{maximal path in} G_\mytimes & \quad \Leftrightarrow \quad P \mytext{maximal path in} G.
\end{align*}
In particular, $(G_\mytimes)^h = (G^h)_\mytimes$; we can thus simply write $G_\mytimes^h$. 
\end{corollary}
\BP The first proposition is due to $G_\mytimes \subseteq G$. Let $P = v_0, \ldots, v_I$ be a maximal path in $G$. According to Corollary \ref{coro_maxPath}, $v_{i-1}, v_i$ is the only path in G from $v_{i-1}$ to $v_i$, and $v_{i-1} v_i \in  A(G_\mytimes)$ follows with \eqref{RdG_formel}. $P$ is thus a path in $G_\mytimes$. If $P'$ is a path in $G_\mytimes$ with $P' \supseteq P$, then $P'$ is according to the first statement a path in $G$, hence $P' = P$, and $P$ is maximal in $G_\mytimes$.

Now let $P$ be a maximal path in $G_\mytimes$. Then $P$ is a path in $G$ which is contained in a maximal path $P'$ in $G$. According to what we have already proven, $P'$ is a path in $G_\mytimes$, hence $P' = P$, and $P$ is maximal in $G$.

\EP

For $G \in \fTa$, we define
\begin{align*}
P_\mytimes & \; \equiv \; ( G \vert_P )_\mytimes \quad \mytext{for every path} P \subseteq V(G), \\
\L(G) & \; \equiv \; \mysetdescr{ P_\mytimes }{ P \in \P_G^h }.
\end{align*}
We have $V(P_\mytimes) = P$, $A( P_\mytimes ) = \mysetdescr{ P_{i-1} P_i }{ i \in \ell(P) }$ for every path $P$ in $G$. In particular, if $P \subseteq V(G)$ and $P' \subseteq V(G')$ are paths in $G, G' \in \fTa$ with $\ell(P) = \ell(P')$, then $P_\mytimes \simeq P'_\mytimes$.

\begin{corollary} \label{coro_z}
Let $G \in \fTa$ and $H \in \fTa \cap \fDr$ with $h_G = h_H$ and $\sigma \in \S(G,H)$. Then $\iotaMxi{A(P_\mytimes)}{\sigma} = \fz_{A(P_\mytimes)}$ for every $P \in \P_G^h$, hence 
\begin{equation} \label{summe_pfad_strict}
2 \cdot h_H \; = \; \sum_{i=1}^{h_G} \iotaxi{\sigma}( P_{i-1}, P_i ) \; = \; \mu_{\sigma \vert_P}( P_\mytimes ).
\end{equation} 
In particular, $\S(G,H) \subseteq \J_{G,H}^{\L(G)}(\sigma)$.
\end{corollary}
\BP Let $\sigma \in \S(G,H)$. For $h_G = 0$, there is nothing to show due to $A( G_\mytimes^h) = \emptyset$ and $A(P_\mytimes) = \emptyset$ for every $P \in \P_G^h$. Let $h_G \geq 1$ and $P \in \P_G^h$. The image $\sigma[ P]$ is a path of length $h_G = h_H$ in $H$, hence a maximal path in $H$, and $\iotaMxi{A(P_\mytimes)}{\sigma} = \fz_{A(P_\mytimes)}$ follows with \eqref{maxPath_iotaFunktion}.

For $\tau \in \S(G,H)$, we have thus $\iotaMxi{A(L)}{\tau} = \fz_{A(L)} = \iotaMxi{A(L)}{\sigma}$ for all $L \in \L(G)$, and $\tau \in \J_{G,H}^{\L(G)}(\sigma)$ follows.

\EP

In our three proofs of implications of type \eqref{zielimplikation} in Section \ref{sec_implications}, the crucial step will be to show that clearly more than $\S(G,H) \subseteq \J_{G,H}^{\L(G)}(\sigma)$ can be said if $G$ and $H$ belong to appropriate classes (Lemma \ref{lemma_R_R_mytimes}, Corollary \ref{coro_SGR_JLGR}, and Theorem \ref{theo_phiG}).

\subsection{The class $\bf \fR$} \label{subsec_R}

In our proofs of implications of type \eqref{zielimplikation} in the next section, we will take $R$ and $S$ from the following class $\fR$ of digraphs:

\begin{definition} \label{def_classR}
We define $\fR$ as the class of digraphs $R \in \fTa \cap \fDr$ with
\begin{equation*}
\M^{\L(R)}(R,R) \cap \S(R,R) \; \not= \; \emptyset.
\end{equation*}
\end{definition}

The following lemma will be the transmission belt putting Theorem \ref{theo_IGR_IGS} into action in Section \ref{sec_implications}. In the terms of Theorem \ref{theo_selExists}, it guarantees for $R \in \fD$ the existence of a strict homomorphism with a selecting arc weight of $G$ within $\H(G,R)$ for every $G \in \fTa$ with $\S(G,R) \not= \emptyset$:

\begin{lemma} \label{lemma_classR}
$R \in \fTa \cap \fDr$ is an element of $\fR$ iff for all $G \in \fTa$ with $h_G = h_R$,
\begin{align*}
\S(G,R) & \not= \emptyset \\
\Rightarrow \quad \M^{\L(G)}(G,R) \cap \S(G,R) & \not= \emptyset.
\end{align*}
\end{lemma}
\BP Let $R \in \fR$ and choose a strict $\rho \in \M^{\L(R)}(R,R)$. Let $G \in \fTa$ with $h_G = h_R$, $P \in \P_G^h$, $\sigma \in \S(G,R)$. Then $\sigma[ P ]$ is a path of length $h_R$ in $R$ and $\rho \circ \sigma$ is strict. Because of
\begin{equation*}
A( \sigma[P]_\mytimes ) \; = \; \mysetdescr{ \sigma( P_{i-1} ) \sigma(P_i) }{ i \in \myNk{h_R} },
\end{equation*}
we have
\begin{align*}
\mu_{\rho \circ \sigma \vert_P}( P_\mytimes )
& \; = \;
\mu_{\rho \vert_{\sigma[P]}}(\sigma[P]_\mytimes) 
\; = \;
\hat{\mu}(\sigma[P]_\mytimes, R ),
\end{align*}
and $\sigma[P]_\mytimes \simeq P_\mytimes$ delivers $\mu_{\rho \circ \sigma \vert_P}(P_\mytimes) = \hat{\mu}(P_\mytimes, R )$. The inverse implication is trivial.

\EP

\begin{proposition} \label{prop_char_R}
A digraph $R \in \fTa \cap \fDr$ is an element of $\fR$ iff 
\begin{align} \label{summenbedingung}
\sum_{i=1}^{\ell(P)} \iota( P_{i-1}, P_i )_R & \; \leq \;
h_R + \ell(P)
\end{align}
for every path $P \subseteq V(R)$. In particular, every flat poset belongs to $\fR$.
\end{proposition}
\BP ``$\Rightarrow$'': Assume that $P \subseteq V(R)$ is a path in $R \in \fTa \cap \fDr$ violating \eqref{summenbedingung}. With $\ell \equiv \ell(P)$ and $n \equiv h_R$, we must have $2 \leq \ell < n$ (use \eqref{summe_pfad_strict} for the last inequality). Let $P' \in \P_R^h$. The homomorphism $\xi : P'_\mytimes \rightarrow R$
\begin{align*}
\xi(P'_i) \; \equiv \; 
\begin{cases}
P_i, & \mytext{if } 0 \leq i \leq \ell; \\
P_\ell, & \mytext{if } \ell < i \leq n
\end{cases}
\end{align*}
is not strict with
\begin{align*}
\mu_\xi(P'_\mytimes)
& \; = \; 
\sum_{i = 1}^n \iota_\xi( P'_{i-1}, P'_i )_R
\; = \; 
n - \ell + \sum_{i = 1}^\ell \iota(  P_{i-1}, P_i )_R
\; > \; 2 \cdot n.
\end{align*}
Because of \eqref{summe_pfad_strict}, the set $ \M( P_\mytimes, R)$ cannot contain a strict homomorphism, thus $R \notin \fR$.

``$\Leftarrow$'':  Assume that \eqref{summenbedingung} holds for every path in $R$. Let $P \in \P_R^h$, $n \equiv h_R$, $\xi \in \H(P_\mytimes, R)$, and define $J$ as the set of the indices $i \in \myNk{n}$ with $\xi(P_{i-1}) = \xi(P_i)$. Then $P' \equiv \mysetdescr{ \xi( P_i )}{ i \in \myNk{n}_0 \setminus J }$ is a path in $R$. Therefore,
\begin{align*}
\mu_\xi(P_\mytimes) 
& \; \; = \;
\sum_{i \in \myNk{n}} \iota( \xi( P_{i-1} ), \xi( P_i ))_R \\
& \; \; = \; \;
\# J + \sum_{i \in \myNk{n} \setminus J} \iota( \xi( P_{i-1} ), \xi( P_i ))_R \\
& \; \stackrel{\eqref{summenbedingung}}{\leq} \;
\# J + n + \ell(P') \; = \; n + \ell(P) \\
& \; \; \leq \; \; 2 \cdot n
\; \stackrel{\eqref{summe_pfad_strict}}{=} \;
\mu_{\id_R \vert_P}(P_\mytimes).
\end{align*}
We conclude $\id_R \vert_P \in \M( P_\mytimes, R)$, and $R \in \fR$ is shown. The addendum about flat posets is a direct consequence.

\EP

\begin{lemma} \label{lemma_C_gehoert_zu_R}
Let $R \in \fTa \cap \fDr$. If the interval $[v, w]_R$ is a path in $R$ for every $vw \in A(R)$, then $R \in \fR$.
\end{lemma} 
\BP Let $P$ be a path in $R$. In the case of $\ell(P) = 0$, $\sum_{i \in \myNk{\ell(P)}} \iota( P_{i-1}, P_i)_R = 0$. Assume $\ell(P) \geq 1$. Then
\begin{equation*}
P' \; \equiv \; \bigcup_{i=1}^{\ell(P)} [ P_{i-1}, P_i ]_R
\end{equation*}
is the concatenation of the paths $[ P_{i-1}, P_i ]_R$, $i \in \myNk{\ell(P)}$, and thus a path in $R$. We get
\begin{align*}
\sum_{i=1}^{\ell(P)} \iota( P_{i-1}, P_i)_R
& \; = \; 
\sum_{i=1}^{\ell(P)} \# [ P_{i-1}, P_i ]_R \\
& \; = \; 
\# P' + \ell(P) - 1 \; = \; \ell(P') + \ell(P)
\; \leq \; h_R + \ell(P),
\end{align*}
and Proposition \ref{prop_char_R} delivers $R \in \fR$.

\EP

\section{Implications of type \eqref{zielimplikation}} \label{sec_implications}

In this section, we prove three implications of type \eqref{zielimplikation} by applying Lemma \ref{lemma_classR} and Theorem \ref{theo_IGR_IGS} on combinations of classes of digraphs $G \in \fTa$ on the one side and $R, S \in \fR$ on the other side.

\subsection{Digraphs $ R \in \fR$ with $R_\mytimes = R^*$} \label{subsec_R_R_mytimes}

\begin{lemma} \label{lemma_R_R_mytimes}
Let $R \in \fTa \cap \fDr$ with  $R_\mytimes = R^*$. Then $R \in \fR$, and for all $G \in \fTa$ with $\S(G,R) \not= \emptyset$, we have for every non-empty set $\L$ of subposets of $G$
\begin{align} \label{MLGR_SLR}
\begin{split}
& \; \M^{\L}(G,R) \\ = \; &
\mysetdescr{ \xi \in \H(G,R) }{ \xi \vert_{V(L)} \in \S(L,G) \mytext{ for all } L \in \L }.
\end{split}
\end{align}
In particular,
\begin{align} \label{SGR_MGGR}
\M(G,R) \; = \; 
\mysetdescr{ \xi \in \H(G,R) }{ \iota_{\xi, A(G^*)} = \fz_{A(G^*)} } \; = \; \S(G,R).
\end{align}
\end{lemma}
\BP Let $G \in \fTa$ and $\xi \in \H(G,R)$. Define
\begin{align*}
J \; \equiv \; \mysetdescr{ vw \in A(G^*) }{ \xi(v) = \xi(w) }
\end{align*}
For $vw \in A(G^*) \setminus J$, we have $\xi(v) \xi(w) \in A(R^*) = A(R_\mytimes)$. There exists a maximal path $P$ in $R_\mytimes$ with $\xi(v) \in P$ and $\xi(w) \in P$. All arcs contained in $A( R \vert_P ) = A(P_\mytimes)$ have the form $P_{i-1} P_i$ with $i \in \ell(P)$, and Corollary \ref{coro_z} delivers $\iota_\xi(v,w) = 2$. Therefore,
\begin{align*}
\mu_\xi(G) \; = \; \sum_{vw \in A(G^*)} \iota_\xi( v,w )
& \; = \; 
2 \cdot \# A(G^*) - \# J \; \leq \; 2 \cdot \# A(G^*)
\end{align*}
with equality iff $J = \emptyset$, hence iff $\xi$ is strict, and \eqref{SGR_MGGR} has been proven for $\S(G,R) \not= \emptyset$.

Now take a subgraph $L \subseteq G$. $\S(G,R) \not= \emptyset$ implies $\S(L,R) \not= \emptyset$. Equation \eqref{SGR_MGGR} delivers $\S(L, R) = \M( L, R)$, and \eqref{MLGR_SLR} follows. For a nonempty set $\L$ of subgraphs of $G$, we have thus $\S(G,R) \subseteq \M^\L(G,R)$, and Lemma \ref{lemma_classR} yields $R \in \fR$.

\EP

\begin{theorem} \label{theo_R_R_mytimes}
Let $R, S \in \fTa \cap \fDr$ with $R^* = R_\mytimes$ and $S^* = S_\mytimes$. Then $R \strH S$ with domain $\fTa$ implies
\begin{equation*}
\# \S(G,R) \leq \# \S(G,S) \quad \mytext{for all } G \in \fTa \mytext{with } \S(G,S) \not= \emptyset.
\end{equation*}
\end{theorem}
\BP Let $G \in \fTa $ with $\S(G,S) \not= \emptyset$. For $\S(G,R) = \emptyset$, there is nothing to prove. Assume $\S(G,R) \not= \emptyset$. Then $\I^{\{G\}}(G,R) = \{ \fz_{A(G^*)} \} = \I^{\{G\}}(G,S)$ according to the first equation in \eqref{SGR_MGGR}. The second one delivers $\S(G,R) = \J_{G,R}^{ \{G\} }(\sigma)$ and $\S(G,S) = \J_{G,S}^{ \{G\} }(\sigma)$ for $\sigma \in \S(G,R)$, and $\# \S(G,R) \leq \# \S(G,S)$ follows with Theorem \ref{theo_IGR_IGS}.

\EP

\subsection{The combination $G \in \fTan$, $R, S \in \fR$} \label{subsec_fTan}

\begin{definition} \label{def_Pxn}
For every $n \in \myN_0$, let $\fTan$ be the class of the finite digraphs $G \in \fTa$ of height $n$ in which every vertex belongs to a path of length $n$.
\end{definition}

For $G \in \fTa$, we have $G^h \in \fTahn{h_G}$, and for $G \in \fTan$, we have $G^h = G$ and $G_\mytimes^h = G_\mytimes$. $\fTahn{1}$ is the class of flat digraphs in $\fTa$ without isolated points. 

\begin{lemma} \label{lemma_fTan}
Let $G \in \fTan$, $n \in \myN_0$, and $H \in \fTa \cap \fDr$ with $h_H \in \{ n, n+1 \}$. Then, for every homomorphism $\xi \in \H(G,H)$,
\begin{equation} \label{strict_auf_SPxH}
\xi \vert_{P} \in \S(P_\mytimes,H) \mytext{ for all } P \in \P_G^h \quad \Leftrightarrow \quad \xi \in \S(G,H).
\end{equation}
In particular, in the case of $h_H = n$, we have for all $\sigma \in \S(G,H)$,
\begin{equation} \label{SGxhH_IGHsigma}
\tau \vert_{V(G)^h} \; \in \; \S( G^h, H ) \quad \mytext{for all } \tau \in \J_{G,H}^{\L(G)}(\sigma).
\end{equation}
\end{lemma}

\BP ``$\Leftarrow$'' in \eqref{strict_auf_SPxH} is trivial. For $n = 0$, also ``$\Rightarrow$'' is trivial, because every mapping from $V(G)$ to $V(H)$ is a strict homomorphism. Assume $n \geq 1$ and let $vw \in A(G^*) \setminus \cup_{P \in \P_G^h} A(P_\mytimes)$. There exist paths $P$ and $P'$ of length $n$ in $G$ with $v \in P$ and $w \in P'$, hence $v = P_{g(v)}$ and $w = P'_{g(w)}$ with $g$ as in Lemma \ref{lemma_gFunktion}. We have $w \not= P_{g(v)+1}$ due to $vw \notin A( P_\mytimes)$.

The sequence $P^+ \equiv P_0, \ldots, P_{g(v)}, P'_{g(w)}, \ldots, P'_n$ is a path in $G$, and the length of $P^+$ is $g(v) + 1 + n - g(w)$. Due to the choice of $vw$, the length of $P^+$ must be less than $n$. We conclude $g(v) + 1 < g(w)$.

Let $\xi \in \H(G,H)$ be a homomorphism fulfilling the condition on the left of \eqref{strict_auf_SPxH}. Then $\xi(P_{i-1}) \xi(P_i) \in A(H^*)$ and $\xi(P'_{i-1}) \xi(P'_i) \in A(H^*)$ for all $i \in \myNk{n}$. In the case of $\xi(v) = \xi(w)$,
\begin{equation*}
\xi(P'_0), \ldots, \xi(P'_{g(w)}), \xi( P_{g(v)+1}), \ldots, \xi(P_n)
\end{equation*}
is a path in $H$ with length $g(w) + n - g(v) \leq h_H \leq n+1$, thus $g(w) \leq g(v) + 1$. Therefore, $\xi(v) \xi(w) \in A(H^*)$, and \eqref{strict_auf_SPxH} is shown.

Now let $h_H = n$, $\sigma \in \S(G,H)$, and $\tau \in \J_{G,H}^{\L(G)}(\sigma)$. For every $P \in \P_G^h$, we have $\iotaMxi{A(P_\mytimes)}{\tau} = \iotaMxi{A(P_\mytimes)}{\sigma} = \fz_{A(P_\mytimes)}$ according to Corollary \ref{coro_z}. With \eqref{iocaChar_strict}, we conclude $\tau \vert_{P} \in \S(P_\mytimes,H)$. We have $G^h \in \fTan$ and $\P_{G^h}^h = \P_G^h$. Now \eqref{strict_auf_SPxH} delivers $\tau \vert_{V(G)^h} \in \S( G^h, H )$.

\EP

\begin{corollary} \label{coro_SGR_JLGR}
For $T \in \fR$ and $G \in \fTan$ with $h_G = h_T = n \in \myN_0$, we have $\S(G,T) = \J_{G,T}^{\L(G)}(\rho)$ for every $\rho \in \S(G,T)$.
\end{corollary}
\BP ``$\subseteq$'' is due to Corollary \ref{coro_z}. Let $\tau \in \J_{G,T}^{\L}(\rho)$. Then, for all $L \in \L(G)$, $\iotaMxi{A(L^*)}{\tau} = \fz_{A(L^*)}$ (Corollary \ref{coro_z} again), and $\tau$ is strict on $L$ according to \eqref{iocaChar_strict}. The equivalence \eqref{strict_auf_SPxH} in Lemma \ref{lemma_fTan} yields $\tau \in \S(G,T)$.

\EP

\begin{theorem} \label{theo_fTan}
(a) Let $R, S \in \fR$ with $h_R = h_S = n \geq 1$. Then $R \strH S$ with domain $\fTa$ implies
\begin{equation*}
\# \S(G,R) \leq \# \S(G,S) \quad \mytext{for all } G \in \fTan \mytext{with } \S(G,S) \not= \emptyset,
\end{equation*}
and  $R \strH S$ with domain $\fP$ implies
\begin{equation*}
\# \S(G,R) \leq \# \S(G,S) \quad \mytext{for all } G \in \fTan \cap \fP.
\end{equation*}

(b) Let $R, S \in \fP$ with $h_R = h_S \leq 1$. Then
\begin{equation*}
R \strH S \quad \Leftrightarrow \quad R \strS S.
\end{equation*}
where $\strH$ and $\strS$ both hold with respect to $\fTa$ or with respect to $\fP$.
\end{theorem}
\BP (a) Let $G \in \fTan$ with $\S(G,S) \not= \emptyset$. For $\S(G,R) = \emptyset$, there is nothing to prove. Assume $\S(G,R) \not= \emptyset$. According to Lemma \ref{lemma_classR}, there exist strict homomorphisms $\sigma \in \M^{\L}(G,R)$ and $\tau \in \M^{\L}(G,S)$, and according to Corollary \ref{coro_z}, we have $\iotaMxi{B}{\sigma} = \fz_B = \iotaMxi{B}{\tau}$ with $B \equiv \cup_{L \in \L(G)} A(L^*)$. We have thus $\iotaMxi{B}{\sigma} \in \I_{G,S}^{\L}$, and for domain $\fTa$, Theorem \ref{theo_IGR_IGS} delivers together with Corollary \ref{coro_SGR_JLGR}
\begin{equation*}
\# \S(G,R) \; = \; \# \J_{G,R}^{\L}(\sigma) 
\; \leq \; 
\# \J_{G,S}^{\L}(\sigma) \; = \; \# \J_{G,S}^{\L}(\tau) \; = \; \# \S(G,S).
\end{equation*}
For domain $\fTa \cap \fP$, repeat the proof with $G \in \fTan \cap \fP$; now $\S(G,S) \not= \emptyset$ is guaranteed.

(b) Let $R, S \in \fP$ with $h_R = h_S \leq 1$. ``$\Leftarrow$'' is due to \eqref{RsqGS_then_RsqS}. Assume $R \not\strS S$ with domain $\fD' \in \{ \fTa, \fP \}$. There exists a digraph $G \in \fD'$ with $\# \S(G,R) > \# \S(G,S)$. Due to $\S(G,R) > 0$, we have $h_G \leq h_R \leq 1$. In the case of $h_G = 0$, we have $\# \H(G,R) = \# \S(G,R) > \# \S(G,S) = \# \H(G,S)$, hence $R \not\strH S$.

Let $h_G = 1$, hence $1 = h_R = h_S$. If $\# V(S) < \#V(R)$, then $\# \H(E,S) < \H(E,R)$ (thus $R \not\strH S$), where $E \in \fD'$ is the singleton with loop. Assume $\# V(S) \geq \#V(R)$. With $i$ denoting the number of isolated points in $G$, we have
\begin{align*}
\# \S(G^h,R) & \; = \; (\# V(R))^{-i} \cdot \# \S(G,R) \\
& \; > \; (\# V(S))^{-i} \cdot \# \S(G,S) \; = \; \# \S(G^h,S).
\end{align*}
$G^h$ is thus an element of $\fTahn{1}$ (of $\fTahn{1} \cap \; \fP$ in the case of $\fD' = \fP$) with $\S(G^h,R) > \S(G^h,S)$. According to Lemma \ref{lemma_C_gehoert_zu_R}, we have $R, S \in \fR$, and because $S$ contains a subset isomorphic to $C_1$, we have $\S(G^h,S) \not= \emptyset$. Now (a) delivers $R \not\strH S$ with domain $\fD'$.

\EP

In the proof of Theorem \ref{theo_fTan}(a), Equivalence \eqref{strict_auf_SPxH} plays an important role. In Proposition \ref{prop_char_fTan} in Section \ref{sec_lemma_fTan_invers}, we will show that for $n \geq 1$, the equivalence characterizes the digraphs in $\fTa$ for which the subgraph formed by the non-isolated points belongs to $\fTan$. We achieve thus no progress by replacing ``for all $G \in \fTan$'' by ``for all $G \in \fTa$ with $h_G = n$ fulfilling \eqref{strict_auf_SPxH}'' in the description of the domain of $\strS$ in Theorem \ref{theo_fTan}(a). We postpone the proof because it does not belong to the main line of this article.

\subsection{The combination $G \in \fTanA$, $R, S \in \fCn$} \label{subsec_fTanA}

Assume $G, H \in \fTa$ with $h_H = h_G$. Lemma \ref{lemma_gFunktion} tells us, that for $v \in V(G)^h$, the position of $\xi(v)$ is independently of $\xi \in \H(G,H)$ highly determined by the position of $v$ in $G$ indicated by $g(v) = g(\xi(v))$. This fact was fundamental for the proof of Lemma \ref{lemma_fTan} and thus for Theorem \ref{theo_fTan}.

We can exploit this determination of images of $v \in V(G)^h$ also for digraphs $G \in \fTa \setminus \fTan$, in which every vertex $v \in V(G)^\noth$ is encapsulated by a set $O(v) \subseteq V(G)^h$. Then the image of $v \in V(G)^\noth$ under a homomorphism $\xi \in \H(G,H)$ is highly determined by the image of the capsule $O(v)$, and if the structure of $\xi[ O(v) ]$ and $\zeta[ O(v) ]$ in $H$ is sufficiently similar for all $\xi, \zeta \in \H(G,H)$, we may hope to achieve a result comparable with Theorem \ref{theo_fTan}. That is the approach in this section. Digraphs $G$ providing a suitable encapsulation of the elements of $V(G)^\noth$ are described in Definition \ref{def_fTanA}, and the sufficient similarity of $\xi[ O(v) ]$ and $\zeta[ O(v) ]$ is enforced by taking $H$ from a suitable subclass of posets contained in $\fR$ (Definition \ref{def_Cn}). In Theorem \ref{theo_fTanA}, we prove a result similar to Theorem \ref{theo_fTan}, but with weaker requirements on $G$ and stronger ones on $R$ and $S$.

Our concept of encapsulation becomes visible in the following definition:

\begin{definition} \label{def_shells}
Let $G \in \fTa$ and $v \in V(G)^\noth$.
\begin{itemize}
\item We call a set $B \subseteq V(G)^h$ a {\em bottom shell of $v$} iff 
\begin{itemize}
\item for every $b \in B$, there exists a path $P$ in $G_\mytimes$ with $P_\bot = b$ and $P_\top = v$;
\item for every path $P$ in $G_\mytimes$ with $P_\bot \in V(G)^h$ and $P_\top = v$, we have $P \cap B \not= \emptyset$.
\end{itemize}
\item We call a set $U \subseteq V(G)^h$ an {\em upper shell of $v$} iff 
\begin{itemize}
\item for every $u \in U$, there exists a path $P$ in $G_\mytimes$ with $P_\bot = v$ and $P_\top = u$;
\item for every path $P$ in $G_\mytimes$ with $P_\bot = v$ and $P_\top  \in V(G)^h$, we have $P \cap U \not= \emptyset$.
\end{itemize}
\end{itemize}
\end{definition}

Now we define:

\begin{definition} \label{def_fTanA}
For $G \in \fTa$, define $\Z(G)$ as the set of connectivity components of $G_\mytimes \vert_{V(G)^\noth}$. We define $\fTanA$ for $n \in \myN_0$ as the class of digraphs $G \in \fTa$ of height $n$ for which
\begin{itemize}
\item for every $v \in V(G)^\noth$, there exists a bottom shell $B(v)$ and an upper shell $U(v)$;
\item  for every $Z \in \Z(G)$, there exist $b_Z, u_Z \in V(G)^h$ with
\begin{equation} \label{LzZUz_in_lZuZ}
\cup_{z \in Z} \; B(z) \; \cup \; Z \; \cup \; \cup_{z \in Z} \; U(z) \quad \subseteq \quad [ b_Z, u_Z ]_{\Tr(G)}.
\end{equation}
\end{itemize}
\end{definition}

\begin{figure}
\begin{center}
\includegraphics[trim = 35 710 220 40, clip]{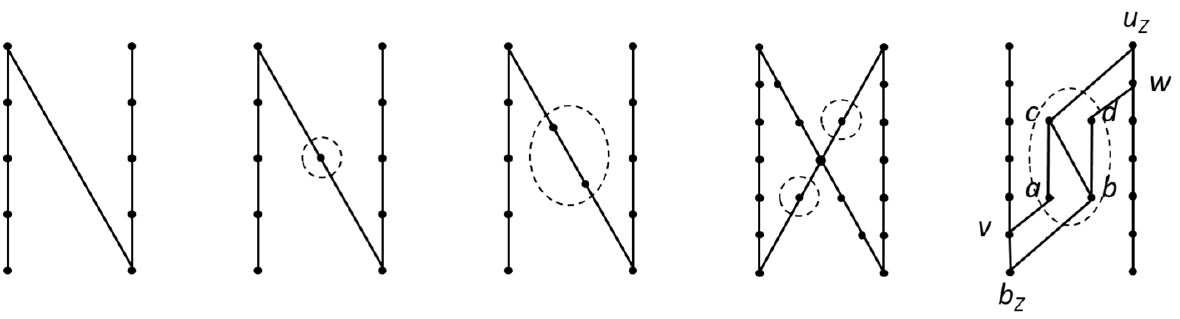}
\caption{\label{fig_fTanA} The Hasse-diagrams of five posets $G \in \fTanA$ and their sets $\Z(G)$. The leftmost poset is an element of $\fTahn{4}$, thus $\Z(G) = \emptyset$; for the other ones, the elements of $\Z(G)$ are encircled.}
\end{center}
\end{figure}

We have $\fTan \subseteq \fTanA$ because of $\Z(G) = \emptyset$ for $G \in \fTan$. Five examples for posets $G \in \fTanA$ are shown in Figure \ref{fig_fTanA}. The leftmost one is an element of $\fTahn{4}$, hence $\Z(G) = \emptyset$. For the second and third one, $\Z(G)$ contains a singleton and the points belonging to the chain of length 1, respectively. For the fourth poset, $\Z(G)$ consists of two singletons, because there is no line connecting these two vertices in $G_\mytimes \vert_{V(G)^\noth}$. For the last poset, $\Z(G)$ contains as single element the vertex set $\{ a, b, c, d \}$.

Let $G \in \fTanA$ and $Z \in \Z(G)$. $b_Z, u_Z \in V(G)^h$ and \eqref{LzZUz_in_lZuZ} imply $B(z) \not= \emptyset$ and $U(z) \not= \emptyset$ for all $z \in Z$. In the most simple case, $B(z) = \{ b_Z \}$ and $U(z) = \{ u_Z \}$ for all $z \in Z \in \Z(G)$, as in the second, third, and fourth poset in Figure \ref{fig_fTanA}. For the fifth poset, $b_Z$ and $u_Z$ are indicated, and we have $B(c) = \{ b_Z, v \}$, $U(b) = \{ u_Z, w \}$, and we can choose $B(a) = \{ v \}$, $U(d) = \{ w \}$.

\begin{definition} \label{def_Cn}
For every $n \in \myN_0$, let $\fCn$ be the class of the posets $R$ in which every maximal path has length $n$ and in which $[ v, w ]_R$ is a path for every $vw \in A(R)$.
\end{definition}
According to Lemma \ref{lemma_C_gehoert_zu_R}, $\fCn \subseteq \fR$. For $R \in \fCn$ and $G \in \fTa$ with $h_G \leq n$, the strict homomorphisms $\lambda_{G,n}$ and $\hat{\lambda}_{G,n}$ defined at the end of Section \ref{subsec_Paths} show $\S(G,R) \not= \emptyset$. As usual for posets, we write $v \leq_Q w$ instead of $vw \in A(Q)$ and $v <_Q w$ instead of $vw \in A(Q^*)$ for a poset $Q$.

\begin{theorem} \label{theo_phiG}
Let $G \in \fTanA$ for $n \in \myN_0$. There exists a number $0 < \phi_G \leq 1$ with
\begin{equation*}
\# \S(G,H) \quad = \quad \phi_G \cdot \# \J_{G,H}^{\L(G)}(\sigma)
\end{equation*}
for every $H \in \fCn$ and every $\sigma \in \S(G,H)$.
\end{theorem}
\BP For $n = 0$, we can use $\phi_G = 1$. Let $n \geq 1$, $H \in \fCn$ and $\sigma \in \S(G,H)$. According to Corollary \ref{coro_z}, we have $\iotaMxi{A(L)}{\sigma} = \fz_{A(L)}$ for every $L \in \L \equiv \L(G)$. The proof is divided in three parts, and $\xi \in \J_{G,H}^\L(\sigma)$ is fixed in the first two of them:
\begin{enumerate}
\item We analyze how $\xi$ maps the interval $[b_Z, u_Z]_G$ to $H$ for $Z \in \Z(G)$.
\item We show that $\xi \vert_{V(G)^\noth}$ can uniquely be represented by a combination of homomorphisms from $G \vert_Z$ to $H$, $Z \in \Z(G)$.
\item We show that this combinatorial representation holds also for strict homomorphisms and derive a formula for $\phi_G$.
\end{enumerate}

1) Let $Z \in \Z(G)$ be selected. We define $\fB \equiv \{ b_Z \} \cup \; \cup_{z \in Z} \; B(z)$, $\fU \equiv \{ u_Z \} \cup \; \cup_{z \in Z} \; U(z)$.

Let $v \in \fB \cup \fU$. The point $v$ belongs to $G^h$; there exist thus a path $P \in \P_G^h$ with $v \in P$, hence $v = P_{g(v)}$ according to Lemma \ref{lemma_gFunktion}. Because of $\iotaMxi{A(P_\mytimes)}{\xi} = \fz_{A(P_\mytimes)}$, the path $\xi[P]$ has length $n$ in $H$ with $\xi(v) = \xi[P]_{g(v)}$. 

Because of $H \in \fCn$ and \eqref{LzZUz_in_lZuZ}, the interval $I(\xi,Z) \equiv [ \xi( b_Z ), \xi( u_Z ) ]_H$ is a path containing $\xi[\fB]$, $\xi[ Z ] $, and $\xi[\fU]$. There exists a maximal path in $H$ containing $I(\xi,Z)$. Owing to $H \in \fCn$, every maximal path $P' \supseteq I(\xi,Z)$ has length $n$, and due to $\xi(v) = \xi[P]_{g(v)}$, Lemma \ref{lemma_gFunktion} yields $\xi(v) = P'_{g(v)}$ for every such path $P'$. The path $I(\xi,Z)$ depends of course on $\xi$, $\sigma$, and $H$, but with
\begin{equation*}
f_Z(w) \; \equiv \; g(w) - g( b_Z ) \quad \mytext{ for all } w \in [ b_Z, u_Z ]_G \cap V(G)^h,
\end{equation*}
we have $\xi(v) = I(\xi,Z)_{f_Z(v)}$; the position of $\xi(v)$ within $I(\xi,Z)$ does thus not depend on $\xi$, $\sigma$, or $H$. Furthermore, with $k_Z \equiv f_Z(u_Z)$, we have $Q(\xi,Z) \equiv H \vert_{I(\xi,Z)} \simeq C_{k_Z}$, and this chain does not depend on $\xi$, $\sigma$, or $H$, too.

2) Let $Z \in \Z(G)$. For every $z \in Z$, we define
\begin{align*}
m_Z(z) & \; \equiv \; \max \mysetdescr{ f_Z(a) }{ a \in B(z) }, \\
M_Z(z) & \; \equiv \; \min \mysetdescr{ f_Z(b) }{ b \in U(z) }.
\end{align*}
Furthermore, for every chain $Q'$ of length $k_Z$, we define $\H_{m_Z,M_Z}(Z, Q')$ as the set of homomorphisms $\theta \in \H(Z, Q')$ for which, for all $z \in Z$,
\begin{equation*}
Q'_{m_Z(z)} \; \leq_{Q'} \; \theta(z) \; \leq_{Q'} \; Q'_{M_Z(z)}.
\end{equation*}
(Here, as in the rest of this proof, we identify a chain with the path formed by its vertices.)

For every $z \in Z$, $b \in B(z)$, and $c \in U(z)$, there exist paths in $G_\mytimes$ from $b$ to $z$ and from $z$ to $c$. Therefore, $\xi \vert_Z \in \H_{m_Z,M_Z} \left( Z, Q(\xi,Z) \right)$ for all $Z \in \Z(G)$. We can thus regard $\xi \vert_{V(G)^\noth}$ as an uniquely determined element of
\begin{equation*}
\prod_{Z \in \Z(G)} \H_{m_Z,M_Z} \left( Z, Q(\xi,Z) \right).
\end{equation*}
Now we show, that every $\Theta$ contained in this Cartesian product 
defines a $\zeta \in \J_{G,H}^\L(\sigma)$ with $\zeta \vert_{V(G)^h} = \xi \vert_{V(G)^h}$ via
\begin{align*}
\zeta(v) \; \equiv \; \begin{cases}
\xi(v), & \mytext{if } v \in V(G)^h; \\
\Theta_Z(v), & \mytext{if } v \in Z \in \Z(G).
\end{cases}
\end{align*}
That is proven if we have shown that for $Z \in \Z(G)$, $\theta \in \H_{m_Z,M_Z} \left( Z, Q(\xi,Z) \right)$, the mapping $\zeta : V(G) \rightarrow V(H)$ defined by
\begin{align*}
\zeta(v) \; \equiv \; \begin{cases}
\xi(v), & \mytext{if } v \in V(G) \setminus Z, \\
\theta(v), & \mytext{if } v \in Z,
\end{cases}
\end{align*}
belongs to $\J_{G,H}^\L(\sigma)$. We start with the proof of $\zeta \in \H(G,H)$. Let $v w \in A(G^*)$. There exists a path $P$ in $G_\mytimes$ from $v$ to $w$, and due to Corollary \ref{coro_pathsInGmytimes}, $P$ is a path in $G$, too.

In the cases $v, w \in V(G) \setminus Z$ and $v, w \in Z$, $\zeta( v ) \leq_H \zeta(w)$ is trivial. Let $v \in V(G) \setminus Z$ and $w \in Z$. In the case of $v \in V(G)^h$, the bottom shell of $w$ delivers an $a \in B(w) \cap P$. Because $H$ is transitive, we get
\begin{align} \label{ungl_kette_1}
\begin{split}
\zeta(v) \; = \; \xi(v) \leq_H \xi(a) 
& \; = \; \; I(\xi,Z)_{f_Z(a)} \\
& \; \leq_H I(\xi,Z)_{m_Z(w)}
\; \stackrel{(*)}{\leq_H} \; \theta(w)
\; = \; \zeta(w). 
\end{split}
\end{align}

Now let $v \in V(G)^\noth \setminus Z$. There exists a $Z' \in \Z(G)$ with $v \in Z'$. Because of $Z' \not= Z$, the sets $Z'$ and $Z$ are different connectivity components of $G_\mytimes \vert_{V(G)^\noth}$, and because $P$ is a path in $G_\mytimes$, there exists an index $i \in \myNk{\ell(P)-1}$ with $P_i \in V(G)^h$. Now the bottom shell of $w$ delivers an index $j$ with $i \leq j \leq \ell(P)-1$ and $P_j \in B(w)$, and as above, we get
\begin{align} \label{ungl_kette_2}
\begin{split}
\zeta(v) \; = \; \xi(v) \leq_H \xi(P_j) 
& \; = \; \; I(\xi,Z)_{f_Z(P_j)} \\
& \; \leq_H I(\xi,Z)_{m_Z(w)}
\; \stackrel{(*)}{\leq_H} \; \theta(w)
\; = \; \zeta(w). 
\end{split}
\end{align}
The proof for $v \in Z, w \in V(G) \setminus Z$ runs analogous,  and $\zeta \in \H(G,H)$ is shown. Now $ \zeta \vert_{V(G)^h}= \xi \vert_{V(G)^h}$ delivers $\iotaMxi{A(L)}{\zeta} = \iotaMxi{A(L)}{\xi} = \fz_{A(L)}$ for all $L \in \L$, hence $\zeta \in \J_{G,H}^\L(\sigma)$.

3) For $\xi \in \J_{G,H}^\L(\sigma)$, we define
\begin{equation*}
\langle \xi \rangle \quad \equiv \quad 
\mysetdescr{ \zeta \in \J_{G,H}^\L(\sigma) }{ \zeta \vert_{V(G)^h} = \xi \vert_{V(G)^h} }.
\end{equation*}
The set $\mysetdescr{ \langle \xi \rangle }{ \xi \in \J_{G,H}^\L(\sigma) }$ is a partition of $\J_{G,H}^\L(\sigma)$. Let $\Xi$ be a representation system. According to the result of part (2) of our proof, we have
\begin{equation*}
\# \langle \xi \rangle \quad = \quad \prod_{Z \in \Z(G)} \# \H_{m_Z,M_Z} \left(Z, \myNk{k_Z}_0 \right),
\end{equation*}
a number not depending on $\xi$, $\sigma$, or $H$.

A necessary condition for the strictness of $\zeta \in \langle \xi \rangle $ is that $\zeta \vert_Z$ is strict for all $Z \in \Z(G)$ with $I(\xi,Z)_{m_Z(z)}  < \zeta(z) < I(\xi,Z)_{M_Z(z)} $ for all $z \in Z$. This condition is sufficient, too. Let $vw \in A(G^*)$. In the case of $v, w \in V(G)^h$, $\zeta(v) <_H \zeta(w)$ follows with \eqref{SGxhH_IGHsigma}, and in the case of $v, w \in Z \in \Z(G)$, $\zeta(v) <_H \zeta(w)$ is due to the strictness of $\zeta \vert_Z$. For $v \in V(G)^h$, $w \in V(G)^\noth$, we get ``$<_H$'' in the inequality marked with $(*)$ in \eqref{ungl_kette_1}; $v \in V(G)^\noth$, $w \in V(G)^h$ is analogous. Finally, for $v \in Z', w \in Z$ with $Z', Z \in \Z(G)$ and $Z' \not= Z$, we have ``$<_H$'' in the inequality marked with $(*)$ in \eqref{ungl_kette_2}.

For every $Z \in \Z(G)$ and for every chain $Q'$ of length $k_Z$, let $\S_{m_Z,M_Z}( Z, Q')$ be the set of homomorphisms $\theta \in \S(Z, Q')$, for which, for all $z \in Z$,
\begin{equation*}
Q'_{m_Z(z)} \; <_{Q'} \; \theta(z) \; <_{Q'} \; Q'_{M_Z(z)}.
\end{equation*}

Due to our last result, there are $\prod_{Z \in \Z(G)} \# \S_{m_Z,M_Z} \left( Z, \myNk{k_Z}_0 \right) $ strict homomorphisms contained in $\langle \xi \rangle $, and also this number does not depend on $\xi$, $\sigma$, or $H$. Because of $\sigma \in \langle \sigma \rangle$, the number is greater than zero.

According to Corollary \ref{coro_z}, we have $\S(G,H) \subseteq \J_{G,H}^\L(\sigma)$. Therefore,
\begin{align*}
\# \J_{G,H}^\L(\sigma) \; & = \; \# \Xi \cdot \prod_{Z \in \Z(G)} \# \H_{m_Z,M_Z} \left( Z, \myNk{k_Z}_0 \right), \\
\# \S(G,H) \; & = \; \# \Xi \cdot \prod_{Z \in \Z(G)} \# \S_{m_Z,M_Z}(Z, \myNk{k_Z}_0), \\
\mytext{hence} \quad \# \S(G,H) \; & = \; 
\# \J_{G,H}^\L(\sigma) \cdot \prod_{Z \in \Z(G)} \frac{\# \S_{m_Z,M_Z} \left( Z, \myNk{k_Z}_0 \right)}{\# \H_{m_Z,M_Z} \left( Z, \myNk{k_Z}_0 \right) }.
\end{align*}

\EP

\begin{theorem} \label{theo_fTanA}
Let $R, S \in \fCn$ with $n \in \myN_0$. Then
\begin{align*}
& \quad R \strH S \mytext{ with domain } \fTa \\
\Rightarrow & \quad R \strS S \mytext{ with domain } \fTanA.
\end{align*}
The implication holds also if we replace $\fTa$ by $\fP$ and $\fTanA$ by $\fTanA \cap \fP$.
\end{theorem}
\BP For $n = 0$, everything is trivial. Assume $n \geq 1$, and let $G  \in \fTanA$ (or $G \in \fTanA \cap \fP)$, $\L \equiv \L(G)$. We have $\S(G,R) \not= \emptyset \not= \S(G,S)$.

Due to $\fCn \subseteq \fR$ and Lemma \ref{lemma_classR}, there exist strict homomorphisms $\rho \in \M^\L(G,R)$ and $\sigma \in \M^\L(G,S)$, and according to Corollary \ref{coro_z}, we have $\iotaMxi{A(L)}{\rho} = \fz_{A(L)} = \iotaMxi{A(L)}{\sigma}$ for all $L \in \L$. In the case of $R \strH S$, Theorem \ref{theo_IGR_IGS} delivers $\# \J_{G,R}^\L(\rho) \leq \# \J_{G,S}^\L(\rho) = \# \J_{G,S}^\L(\sigma)$ for both choices of the domain, and Theorem \ref{theo_phiG} yields
\begin{align*}
\# \S(G,R) \; & = \; \phi_G \cdot \# \J_{G,R}^\L(\rho) 
\; \leq \; 
\phi_G \cdot \# \J_{G,S}^\L(\sigma) \; = \; \# \S(G,S).
\end{align*}

\EP

\section{The inverse of Lemma \ref{lemma_fTan}} \label{sec_lemma_fTan_invers}

In the proof of Theorem \ref{theo_fTan}(a), Equivalence \eqref{strict_auf_SPxH} was the key. We show that for $n \geq 1$, it characterizes the digraphs in $\fTa$ for which the subgraph formed by the non-isolated points belongs to $\fTan$:

\begin{proposition} \label{prop_char_fTan}
Let $G \in \fTa$ be a digraph without isolated points and $h_G \geq 1$. Assume that for every $H \in \fTa \cap \fDr$ with $h_H \in \{ h_G, h_G + 1 \}$ and every $\xi \in \H(G,H)$
\begin{equation*}
\xi \vert_{P} \in \S(P_\mytimes,H) \mytext{ for all } P \in \P_G^h \quad \Leftrightarrow \quad \xi \in \S(G,H).
\end{equation*}
Then $G \in \fTahn{h_G}$.
\end{proposition}
\BP Let $n \equiv h_G$ We define for every $v \in G$ and every $P \in \P_G^h$
\begin{itemize}
\item $K^{in}(v,P)$ is the set of indices $i \in \myNk{n}_0$ for which a path $P'$ exists in $G$ from $P_i$ to $v$,
\item $K^{out}(v,P)$ is the set of indices $i \in \myNk{n}_0$ for which a path $P'$ exists in $G$ from $v$ to $P_i$,
\end{itemize}
and
\begin{align*}
\kappa^{in}(v) & \quad \equiv \quad \max \left\{ \; \{ 0 \} \; \cup \; \bigcup_{P \in \P_G^h} K^{in}(v,P) \right\}, \\
\kappa^{out}(v) & \quad \equiv \quad \min \left\{ \; \{ n \} \; \cup \; \bigcup_{P \in \P_G^h} K^{out}(v,P) \right\}.
\end{align*}
$ \kappa^{in}$ and $\kappa^{out}$ are both elements of $\H(G, \myNk{n}_0)$. Let $P \in \P_G^h$ and $v \in P$. Then $v = P_{g(v)}$, thus $g(v) \leq \kappa^{in}(v)$, where $g$ is the mapping defined in Lemma \ref{lemma_gFunktion}. Assume ``$<$''. There exists a path $P' \in \P_G^h$ for which a path $z_0, \ldots , z_I$ exists from $P'_{\kappa^{in}(v)}$ to $v$. We have $P'_{\kappa^{in}(v)} \not= v$, and
\begin{equation*}
Q \quad \equiv \quad P'_0, \ldots , P'_{\kappa^{in}(v)}, z_1, \ldots , z_I, P_{g(v)+1}, \ldots , P_n
\end{equation*}
is a path in $G$ with $\kappa^{in}(v) + 1 + n - g(v) \leq \ell(Q) \leq n$, hence $\kappa^{in}(v) < g(v)$. Contradiction! Therefore, $\kappa^{in}(v) = g(v)$, thus $\kappa^{in} \vert_P \in \S( P_\mytimes, \myNk{n}_0 )$ for all $P \in \P_G^h$. $\kappa^{out} \vert_P \in \S( P_\mytimes, \myNk{n}_0 )$ for all $P \in \P_G^h$ is shown in the same way.

Now assume $G \notin \fTan$. There exists a $v \in V(G)$ which does not belong to any path in $\P_G^h$. Because $G$ does not contain isolated points, at least one of the sets
\begin{align*}
\Nin(v)^o & \quad \equiv \quad \Nin(v) \setminus \{ v \} , \\
\Nout(v)^o & \quad \equiv \quad
\Nout(v) \setminus \{ v \}
\end{align*}
is non-empty. Assume $\Nin(v)^o \not= \emptyset$. If $P$ is a path in $G$ of length greater/equal 1 ending in $v$, then $P_{\ell(P)-1} \in \Nin(v)^o$, and there is no path $P' \in \P_G^h$ with $ P_{\ell(P)-1}, v \in P'$. Therefore,
\begin{equation*}
\kappa^{in}(v) \quad = \quad \max \mysetdescr{ \kappa^{in}(w) }{ w \in \Nin(v)^o },
\end{equation*}
and $\kappa^{in}$ is not strict. And in the case of $\Nout(v)^o \not= \emptyset$, the homomorphism $\kappa^{out}$ is not strict.

\EP


\begin{thebibliography}{xx}
\bibitem{Aho_etal_1972} A. V. Aho, M. R. Garey, and J. D. Ullman: The transitive reduction of a directed graph. {\em SIAM J. Comput.} {\bf 1} (1972), 131--137.

\bibitem{Atserias_etal_2021} A. Atserias, P. G. Kolaitis, and W.-L. Wu: On the Expressive Power of Homomorphism Counts. arXiv:2101.12733v1 (2021).

\bibitem{Borgs_etal_2006} C. Borgs, J. Chayes, L. Lov\'{a}sz, V. T. S\'{o}s, and K. Vesztergombi: Counting Graph Homomorphisms. In: M. Klazar, J. Kratochvil, M. Loebl, J. Matou\v{s}ek, R. Thomas, and P. Valtr eds., {\em Topics in Discrete Mathematics} (Algorithms and Combinatorics {\bf 26}), Springer (2006), 315--371.

\bibitem{Borgs_etal_2008} C. Borgs, J. Chayes, L. Lov\'{a}sz, V. T. S\'{o}s, and K. Vesztergombi: Convergent sequences of dense graphs I: Subgraph frequencies, metric properties and testing. {\em Advances in Mathematics} {\bf 219} (2008), 1801--1851.

\bibitem{Borgs_etal_2012} C. Borgs, J. Chayes, L. Lov\'{a}sz, V. T. S\'{o}s, and K. Vesztergombi: {\em Convergent sequences of dense graphs II. Multiway Cuts and Statistical Physics}. Ann. of Math. {\bf 176} (2012), 151--219. 

\bibitem{Cai_Govorov_2020} J. Cai and A. Govorov: On a theorem of Lov\'{a}sz that $hom(\cdot,H)$ determines the isomorhphism type of $H$. ITCS (2020).

\bibitem{aCampo_2018} F. a Campo: Relations between powers of Dedekind numbers and exponential sums related to them. {\em J. Int. Seq.}  {\bf 21} (2018), Article 18.4.4.

\bibitem{aCampo_toappear_0} F. a Campo: Criteria for the less-equal-relation between partial Lov\'{a}sz-vectors of digraphs. arXiv:2008.03279v2 (2020).

\bibitem{aCampo_toappear_1} F. a Campo: Generalized one-to-one mappings between homomorphism sets of digraphs. arXiv:1906.11758v4 (2020).

\bibitem{Chaudhuri_Vardi_1993} S. Chaudhuri and M. Y. Vardi: Optimization of {\em Real} conjunctive queries. In C. Beeri (ed.): {\em Proceedings of the Twelfth ACM SIGACT-SIGMODSIGART Symposium on Principles of Database Systems, May 25-28, 1993, Washington, DC} (1993), 59–70.

\bibitem{Dell_etal_2018} H. Dell, M. Grohe, and G. Rattan: Lov\'{a}sz meets Weisfeiler and Leman. In I. Chatzigiannakis, C. Kaklamanis, D. Marx, and D. Sannella (eds.): {\em 45th International Colloquium on Automata, Languages, and Programming, ICALP 2018, July 9-13, 2018, Prague}, LIPIcs {\bf 107} (2018), 40:1–40:14.

\bibitem{Dvorak_2010} Z. Dvo\v{r}\'{a}k: On Recognizing Graphs by Numbers of Homomorphisms. {\em Journal of Graph Theory} {\bf 64} (2010), 330--342.

\bibitem{Erne_1982} M. Ern\'{e}: {\em Einf\"{u}hrung in die Ordnungstheorie}. Bibliographisches Institut, Mannheim, Wien, Z\"{u}rich 1982.

\bibitem{Fisk_1995} S. Fisk: Distinguishing graphs by the number of homomorphisms. {\em Discussiones
Mathematicae - graph theory} {\bf 15} (1995), 73--75.

\bibitem{Freedman_etal_2007} M. Freedman, L. Lov\'{a}sz, and A. Schrijver: Reflection positivity, rank connectivity, and homomorphism of graphs. {\em J. Amer. Math. Soc.} {\bf 20} (2007), 37--51.

\bibitem{Hell_Nesetril_2004} P. Hell and J. Ne\v{s}et\v{r}il: {\em Graphs and Homomorphisms.} Oxford Lecture Series in Mathematics and its Applications {\bf 28}, Oxford University Press (2004).

\bibitem{Lovasz_1967} L. Lov\'{a}sz: Operations with structures. {\em Acta Math. Acad. Sci. Hungar.}  {\bf 18} (1967), 321--328.

\bibitem{Lovasz_2006} L. Lov\'{a}sz: The rank of connection matrices and the dimension of graph algebras. {\em European Journal of Combinatorics} {\bf 27} (2006), 962--970.

\bibitem{Lovasz_Szegedy_2008} L. Lov\'{a}sz and B. Szegedy. Contractors and connectors of graph algebras. {\em Journal of Graph Theory} {\bf 60} (2009), 11--30.

\bibitem{Schrijver_2009} A. Schrijver: Graph invariants in the spin model. {\em J. Combin. Theory B} {\bf 99} (2009), 502--511.

\bibitem{Schroeder_2016} B. Schr\"{o}der: {\em Ordered Sets. An Introduction with Connections from Combinatorics to Topology.} Birkh\"{a}user (2016).

\bibitem{Topsoe_1974} F. Tops{\o}e: {\em Informationstheorie}. B. G. Teubner, Stuttgart 1974.

\end{thebibliography}
\end{document}